\documentclass[11pt]{article}

\usepackage{amsmath,amsfonts,amssymb,amsthm}
\usepackage[margin=1in]{geometry}
\usepackage[colorlinks]{hyperref}
\usepackage{graphicx}
\usepackage{algorithm}
\usepackage{algorithmicx}
\usepackage{authblk}
\usepackage{algpseudocode}
\usepackage{bbm}
\usepackage{comment}
\usepackage{makecell}
\usepackage{accents}
\usepackage{adjustbox}
\usepackage{booktabs}
\usepackage{caption}
\usepackage{amsmath}
\usepackage{subcaption} 

\newcommand{\R}{{\mathbb{R}}}

\usepackage{bm}

\theoremstyle{definition}

\newtheorem{definition}{Definition}[section]

\newcommand{\Be}{\bm{e}}

\newcommand{\Bp}{\bm{p}}

\def\R{{\mathbb{R}}}

\def\bP{\mathbb{P}}

\def\bw{\mathbf{w}}
\def\bx{\mathbf{x}}

\def\by{\mathbf{y}}

\newcommand{\quotes}[1]{``#1''}

\DeclareMathOperator*{\argmin}{arg\,min}

\begin{document}


\title{%
A Finite Expression Method for Solving High-Dimensional Committor Problems}
\author[1]{Zezheng Song\thanks{Email: \texttt{zsong001@umd.edu}}}
\author[1]{Maria K. Cameron\thanks{Email: \texttt{mariakc@umd.edu}}}
\author[1]{Haizhao Yang\thanks{Email: \texttt{hzyang@umd.edu}}}
\affil[1]{Department of Mathematics, University of Maryland, College Park, MD 20742, USA}
\date{}
\maketitle

\begin{abstract}
    Transition path theory (TPT) is a mathematical framework for quantifying
    rare transition events between a pair of selected metastable states $A$ and $B$. Central to TPT is the committor function, which describes the probability to hit the metastable state $B$ prior to $A$ from any given starting point of the phase space. Once the committor is computed, the transition channels and the transition rate can be readily found. The committor is the solution to the backward Kolmogorov equation with appropriate boundary conditions. However, solving it is a challenging task in high dimensions due to the need to mesh a whole region of the ambient space. In this work, we explore the finite expression method (FEX, Liang and Yang (2022)) as a tool for computing the committor. FEX approximates the committor by an algebraic expression involving a fixed finite number of nonlinear functions and binary arithmetic operations. The optimal nonlinear functions, the binary operations, and the numerical coefficients in the expression template are found via reinforcement learning. The FEX-based committor solver is tested on several high-dimensional benchmark problems. It gives comparable or better results than neural network-based solvers. Most importantly, FEX is capable of correctly identifying the algebraic structure of the solution which allows one to reduce the committor problem to a low-dimensional one and find the committor with any desired accuracy.
\end{abstract}

\textbf{Keywords}  Committor Functions; Rare Events; Finite Expression Method; High Dimensions; Deep Neural Network; Symbolic Learning.

\textbf{AMS} 65N99; 68T07.

\section{Introduction}
Understanding the transition events of a stochastic system between disjoint metastable states is of great importance in many branches of applied sciences~\cite{zahn2004nucleation,zhao2015phase,grant2010large,berteotti2009protein,okuyama1998transition}. Examples of such transition events include conformational changes of biomolecules and dislocation dynamics in crystalline solids. In this work, we focus on the overdamped Langevin process as the underlying dynamics
\begin{equation*}
    \label{eqn:overdamped_langevin}
    d \bx_t=-\nabla V\left(\bx_t\right) d t+\sqrt{2 \beta^{-1}} d \bw_t,
\end{equation*}
where $\bx_t \in \Omega \subset \mathbb{R}^d$ is the state of the system at time $t$, $V: \R^d \rightarrow \R$ is a smooth and coercive potential function, 
$\beta^{-1} = k_{B}T$ is the absolute temperature times the Boltzmann constant, and $\mathbf{w}_t$ is the standard $d$-dimensional Brownian motion. 
The invariant probability density for the overdamped Langevin dynamics~\eqref{eqn:overdamped_langevin} is given by
\begin{equation}
\label{eqn:inv_pdf}
    \rho(\bx) = Z_{\beta}^{-1}\exp(-\beta V(\bx)),\quad Z_{\beta} =\int_{\mathbb{R}^d}\exp(-\beta V(\bx))d\bx.
\end{equation}
In real-world applications, the dimension $d$ is typically high, leading to the major numerical difficulty of concern in this paper.

Transition path theory (TPT)~\cite{vanden2006towards,vanden2010transition} is a mathematical framework for the quantitative description of statistics of transition events. The committor function is a central object in TPT. The reaction rate, the reaction current, i.e., the vector field delineating the transition process, and the density of transition paths are expressed in terms of the committor. For two 
disjoint regions $A$ and $B$ in $\Omega$ chosen by the user, the committor function $q(\bx)$ is defined as 
\begin{equation*}
\label{eqn:committor_func}
    q(\bx) = \bP(\tau_B < \tau_A | \bx_0 = \bx),
\end{equation*}
where $\tau_A$ and $\tau_B$ are the hitting times for the sets $A$ and $B$, respectively. 
The committor function is the solution to the boundary-value problem (BVP) for the backward Kolmogorov equation
\begin{align}   
    &\beta^{-1} \Delta q -\nabla V \cdot \nabla q=0,\quad \bx\in\Omega_{AB}:=\Omega \backslash(\bar{A} \cup \bar{B}), \label{eqn:backward_kolmogorov}\\
    &\left.q(\bx)\right|_{\partial A}=0,\quad\left.q(\bx)\right|_{\partial B}=1,\notag\\
    &\frac{\partial q}{\partial\hat{n}} = 0,\quad \bx\in\partial\Omega,\notag
\end{align}
where $\hat{n}$ is the outer unit normal vector.

An analytical solution to equation~\eqref{eqn:backward_kolmogorov} can be found only in special cases. Otherwise, BVP \eqref{eqn:backward_kolmogorov} must be solved numerically. However, the curse of dimensionality makes traditional numerical schemes, such as finite difference and finite element methods, prohibitively expensive when $d > 3$. To address this issue, alternative methods based on ideas borrowed from data science and machine learning have been emerging. 


\subsection{Prior Approaches to Solving the Committor Problem in High Dimensions}
\label{sec:priorwork}
 To the best of our knowledge, the first high-dimensional committor solver was introduced by Lai and Lu~\cite{lai2018point} relying on the assumption that the dynamics were concentrated near a low-dimensional manifold. The authors constructed a linear system for the committor problem via introducing local meshes on patches of the point cloud and corrected the entries of the resulting stiffness matrix afterward to make it symmetric.  While this idea is elegant, we found that it suffers from a lack of robustness when the intrinsic dimension of the system varies throughout the phase space.

Neural network-based committor solvers were introduced in~\cite{khoo2019solving,lilinren2019,rotskoff2020learning}. They exploit the fact that BVP \eqref{eqn:backward_kolmogorov} admits a variational formulation
\begin{equation}
    \label{eqn:variational_form}
   \argmin_{\substack{f\in\mathcal{C}^1(\Omega_{AB})\\
   \left.f(\bx)\right|_{\partial A}=0,~\left.f(\bx)\right|_{\partial B}=1}} \int_{\Omega_{AB}}|\nabla f(\bx)|^2 \rho(\bx) d \bx,
\end{equation}
where $\rho$ is the invariant density \eqref{eqn:inv_pdf}.
The committor is approximated by a solution model involving a neural network which is trained to minimize the objective function \eqref{eqn:variational_form}. In \cite{lilinren2019}, the boundary conditions for the committor are built in the solution model, while in \cite{khoo2019solving} they are enforced by means of penalty functions. The solution models in these works are also different. The model \cite{khoo2019solving} involves Green's functions to facilitate accurate approximation of the committor at high temperatures, while the model in \cite{lilinren2019} is advantageous for low temperatures as it allows for training points generated using enhanced sampled algorithms. The accuracy of both of these solvers is limited by the accuracy of Monte Carlo integration as the integral in \eqref{eqn:variational_form} is approximated as a Monte Carlo sum.

Gao et al.~\cite{gao2023transition} presented a data-driven approach to efficiently compute transition paths using optimal control theory and discrete Markov processes. Their method
excels in identifying low-dimensional sturctures, enhancing the simulation of rare events. Similarly,~\cite{gao2023optimal} introduced a stochastic optimal control formulation to compute transition paths in Markov jump processes over an infinite time horizon. It establishes a framework using the Girsanov transformation to minimize the relative entropy of path measures, effectively identifying the committor function as a critical tool for determining optimal paths. By using diffusion maps to learn reaction coordinates, the authors of~\cite{gao2023data} reduced the problem’s dimensionality and proposed an unconditionally stable finite volume scheme for the corresponding Fokker-Planck equation on the manifold.

A committor solver based on tensor train representations was introduced in~\cite{chen2023committor}. This approach is shown to be suitable even for committor problems for discretized stochastic partial differential equations and yields an accurate solution if the natural computational domain is box-shaped. However, it can suffer from low accuracy when the geometry of the problem is more complicated so that the involved functions (the committor, the invariant density, etc.) cannot be represented accurately by linear combinations of a few basis polynomials.

Finally, techniques based on diffusion maps \cite{banisch2020,evans2021computing,evans2022JCP} are suitable for finding committors in high dimensions provided that the intrinsic dimension of the problem is not very high, e.g. $d=4$. These approaches feature approximations to the backward Kolmogorov operator in \eqref{eqn:backward_kolmogorov} using diffusion kernels and various renormalizations accounting for enhanced sampling data and anisotropic and position-dependent diffusion resulting from the use of physically motivated dimensional reduction. The committor function is found on the data points. The accuracy of this approach is limited by Monte Carlo integration and a finite bandwidth parameter inherent to the construction.

\subsection{The Goal and a Brief Summary of Main Results}
Liang and Yang recently introduced a novel approach to solving high-dimensional PDEs named the \emph{finite expression method (FEX)}~\cite{liang2022finite}. This method is capable of finding solutions to PDEs  with machine precision provided that the exact solutions can be approximated as algebraic expressions involving a relatively small number of nonlinear functions and binary operations. This was demonstrated on a series of test problems of dimensions up to 50. In FEX, the problem of solving a PDE numerically is transformed into a mixed optimization problem involving both combinatorial optimizations for operator selection and continuous optimization for trainable parameters.

The goal of this work is to investigate FEX as a committor solver, understand its strengths and limitations, and identify routes for further improvements. 
The committor problem arising in applications is often high-dimensional. Furthermore, the phase space $\Omega$ of the underlying stochastic process is often unbounded, and hence the computational domain needs to be chosen by the user wisely. The sets $A$ and $B$ are often chosen to be balls, ellipsoids or potential energy sublevel sets surrounding two selected local minima of the potential energy. This choice always makes the computational domain nonconvex and the functional dependence of the solution on phase variables nontrivial. 
Therefore, the committor problem in high dimensions presents a challenge for FEX that is worthy of exploring.  We emphasize, that the test problems in \cite{liang2022finite} all had exact solutions given by short formulas -- the fact that is not true for the committor problem.


In this work, we adapt FEX approach to the committor problem, discuss its setup and implementation, and apply it to a number of test problems. The test problems include the double-well potential with hyperplane boundaries as in \cite{khoo2019solving,chen2023committor}, the double-well potential with sublevel sets boundaries, concentric spheres as in \cite{khoo2019solving}, rugged Mueller's potential with as in \cite{khoo2019solving}, and butane.
Our findings are the following.
\begin{enumerate}
\item On the benchmark test problems, FEX performs comparably or better than neural network-based solvers. 
\item Remarkably, FEX is capable of capturing the algebraic structure of the solution, i.e. identifying variables or combinations of variables on which the solution depends and does not depend. This allows for dimensional reduction and the use of traditional highly accurate methods such as Chebyshev spectral methods or finite element methods to find the solution with the desired precision. This ability is unique to FEX. The other high-dimensional committor solvers mentioned in Section \ref{sec:priorwork} do not have it.
\end{enumerate}


The rest of the paper is organized as follows. The FEX algorithm is detailed in Section~\ref{sec:proposed_method}. The applications of FEX to four benchmark test problems are presented in Section~\ref{sec:numericalexp}. The conclusion is formulated in Section~\ref{sec:summary}. Additional numerical details are found in Supplementary Information.

\section{
The Finite Expression Method
}
\label{sec:proposed_method}
\begin{figure}[ht]
    \centering
    \includegraphics[width=0.80\linewidth]{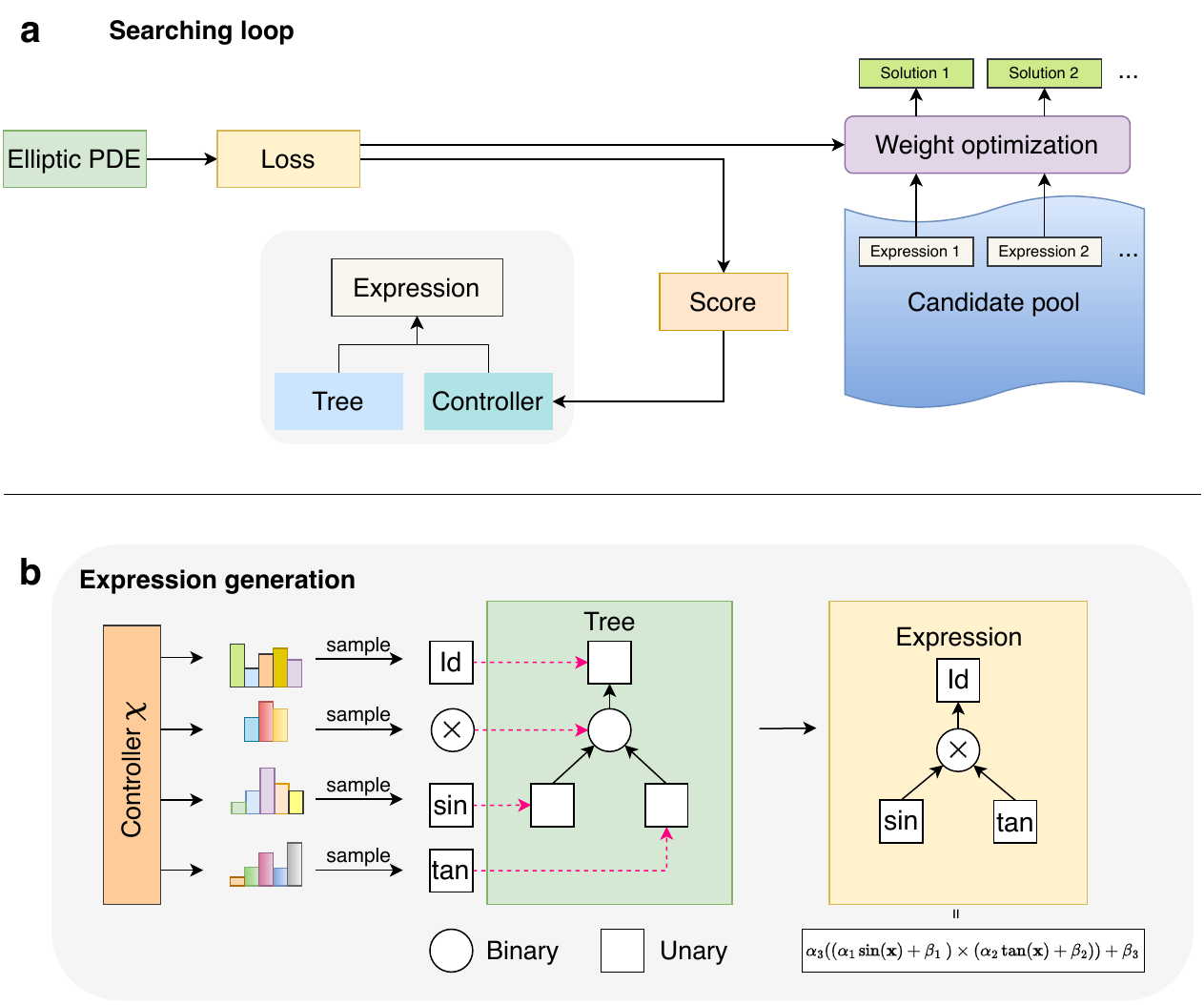}
    \caption{Representation of the components of our FEX implementation. (a) The searching loop for the symbolic solution encompasses multiple stages, namely expression generation, score computation, controller update, and candidate optimization. (b) Illustration of the expression generation with a binary tree and a controller $\bm{\chi}$. The controller produces probability mass functions for each node of the tree, enabling the sampling of node values. Furthermore, we incorporate learnable scaling and bias parameters to generate expressions based on the predefined tree structure and the sampled node values. }
    \label{fig:workflow}
\end{figure}

The finite expression method (FEX)~\cite{liang2022finite} is a brand-new approach to solving PDEs.
It seeks the solution in the form of a mathematical expression with a predefined finite number of operators.  The workflow of FEX is depicted in Figure~\ref{fig:workflow}. To implement FEX, the user needs to choose a binary tree of a finite depth, usually between 4 and 6. This tree is realized as a computer algebra expression. Each tree node is associated with one operator, unary or binary. The input variables are propagated through the tree by passing them through the leaf nodes. Each operator is equipped with a set of parameters defining an affine operator acting on its input argument. Lists of possible unary and binary operators are supplied by the user.
Once the operators and parameters have been assigned to all tree nodes, the binary tree can generate a mathematical expression. The set of operators and parameters minimizing the loss functional measuring how well the expression fits the PDE is sought by solving the mixed combinatorial optimization problem.


To solve this mixed optimization problem, FEX utilizes a combinatorial optimization method to search for the optimal choice of operators and a continuous optimization method to identify the associated parameters of tree nodes. In combinatorial optimization, a reinforcement learning approach is adopted to further reformulate combinatorial optimization into continuous optimization over probability distributions. To accomplish this, FEX introduces a controller network that outputs the probabilities of selecting each operator at each node of the tree. As a result, the problem of selecting the best operators for the tree nodes is reformulated as the problem of identifying the best controller network that can sample the best operators. The optimization for the best controller is a continuous optimization problem. The parameters of the controller network are trained using policy gradient methods to maximize the expected reward in the reinforcement learning terminology, which is equivalent to minimizing the loss function. Therefore, the optimal choice of operators can be determined by sampling from the output probabilities of the best controller network with high probability.


\subsection{The Functional Space of Finite Expressions}
\label{sec:abstractframework}

FEX approximates the solution of a PDE in the space of functions with finitely many operators. Therefore, it is important to formally define the functional space in which the solution is sought.
\begin{definition}[Mathematical expression~\cite{liang2022finite}]\label{def:math_expression}\label{def:expression}
A mathematical expression is a combination of symbols, which is well-formed by syntax and rules and forms a valid function. The symbols include operands (variables and numbers), operators (e.g., \quotes{+}, \quotes{sin}, integral, derivative), brackets, and punctuation.
\end{definition}

\begin{definition}[$k$-finite expression~\cite{liang2022finite}]\label{def:k-finite}
A mathematical expression is called a $k$-finite expression if the number of operators in this expression is $k$.
\end{definition}

\begin{definition}[Finite expression method~\cite{liang2022finite}]\label{def:FEX}
The finite expression method is a methodology to solve a PDE numerically by seeking a finite expression such that the resulting function solves the PDE approximately.
\end{definition}

We denote $\mathbb{S}_k$ the functional space that consists of functions formed by finite expressions with the number of operators less or equal to $k$.

\subsection{The Mixed Combinatorial Optimization Problem in FEX}\label{sec:error} 
The loss functional $\mathcal{L}$ in FEX is problem-dependent.  Reasonable choices include
the least-squares loss as in~\cite{lagaris1998artificial,sirignano2018dgm,nn1994}, a variation formulation as in~\cite{yu2018deep,CiCP-29-1365}, and a 
weak formulation as in~\cite{chen2020friedrichs,zang2020weak}. In this work, we use the fact that the committor problem \eqref{eqn:backward_kolmogorov} admits the variational formulation \eqref{eqn:variational_form} and choose the variation loss functional as in~\cite{khoo2019solving}  given by
\begin{equation}
    \label{eqn:regularized_variational_form}
    \mathcal{L}(u) = \int_{\Omega_{AB}}\left\|\nabla u(\mathbf{x})\right\|^2 \rho(\mathbf{x}) d \mathbf{x}+\tilde{c} \int_{\partial A} u(\mathbf{x})^2 d m_{\partial A}(\mathbf{x})+\tilde{c} \int_{\partial B}\left(u(\mathbf{x})-1\right)^2 d m_{\partial B}(\mathbf{x}).
\end{equation}
Here $\rho$ is the invariant density \eqref{eqn:inv_pdf} and $m_{\partial A}, m_{\partial B}$ are the user-chosen measures on the boundaries $A$ and $B$, respectively. The advantage of this form is that it requires evaluation of only the first derivatives of solution candidates, unlike the least squares loss. This saves the runtime. The boundary conditions in \eqref{eqn:backward_kolmogorov} are enforced by means of penalty terms which is simple and convenient. In FEX, the solution is found by solving the mixed combinatorial optimization problem 
\begin{align}
    \min_{u \in \mathbb{S}}\mathcal{L}(u),
    \label{eqn:co}
\end{align}
where the solution space $\mathbb{S}_{\sf FEX}\subset\mathbb{S}_k$ will be elaborated in Section \ref{sec:binary_tree}.

\subsection{Implementation of FEX}
\label{sec:alg}

The computational workflow of FEX starts with the construction of a binary tree. Each tree node contains either a unary or binary operator. The solution candidate $u$ is obtained by the evaluation of the function represented by the binary tree.  Next, the mixed combinatorial optimization ~\eqref{eqn:co} is applied to adaptively select the optimal operators in all tree nodes. The goal of the mixed combinatorial optimization is to identify operators that can recover the structure of the true solution. FEX  is summarized in Algorithm~\ref{alg:workflow}. Details of its implementation are found in the supplementary materials.

\begin{algorithm}[ht]  
    \caption{Coefficient filtering FEX with a fixed tree to solve committor functions}  
    \label{alg:workflow} 
    \textbf{Input:} PDE and the associated functional $\mathcal{L}$; A tree $\mathcal{T}$; Searching loop iteration $T$; Coarse-tune iteration $T_1$ with Adam; Coarse-tune iteration $T_2$ with BFGS; Fine-tune iteration $T_3$ with Adam; Pool size $K$; Batch size $N$; Coefficient filtering threshold $\tau$.
    
    \textbf{Output:} The solution $u(\bm{x};\mathcal{T},\hat{\Be},\hat{\bm{\theta}})$.
    \begin{algorithmic}[1]
        \State Initialize the agent $\bm{\chi}$ for the tree $\mathcal{T}$
        \State $\mathbb{P} \gets \{\}$
        \For{$\_$ from $1$ to $T$}
        \State Sample $N$ sequences $\{\Be^{(1)}, \Be^{(2)}, \cdots, \Be^{(N)}\}$ from $\bm{\chi}$
        \For{$n$ from $1$ to $N$}
        \State Optimize $\mathcal{L} (u(\bm{x};\mathcal{T},\Be^{(n)},\bm{\theta}))$ by coarse-tune with $T_1+T_2$ iterations.
        \State Compute the reward $R(\Be^{(n)})$ of $\Be^{(n)}$
        \If{$\Be^{(n)}$ belongs to the top-$K$ of $S$}
        \State $\mathbb{P}$.append($\Be^{(n)}$)
        \State $\mathbb{P}$ pops some $\Be$ with the smallest reward when overloading
        \EndIf
        \EndFor
        \State Update $\bm{\chi}$ using~\eqref{eqn:mcmcrisk}
        \EndFor
            \For{$\Be$ in $\mathbb{P}$}  
                \State Fine-tune $\mathcal{L} (u(\bm{x};\mathcal{T},\Be,\bm{\theta}))$ with $T_3$ iterations, apply coefficient filtering with threshold $\tau$.
            \EndFor
        \State \Return the expression with the smallest fine-tune error. 
    \end{algorithmic}
    \textbf{Description:} This algorithm presents a procedure for solving committor functions using the Coefficient Filtering FEX approach. The algorithm utilizes a fixed tree structure and incorporates coefficient filtering to remove trivial coefficients below the threshold $\tau$. By iteratively optimizing the associated functional $\mathcal{L}$, the algorithm identifies and fine-tunes the expression with the smallest fine-tune error, ultimately providing the solution $u(\bm{x};\mathcal{T},\hat{\Be},\hat{\bm{\theta}})$.
\end{algorithm}

\begin{figure}
    \centering
    \includegraphics[width=0.8\linewidth]{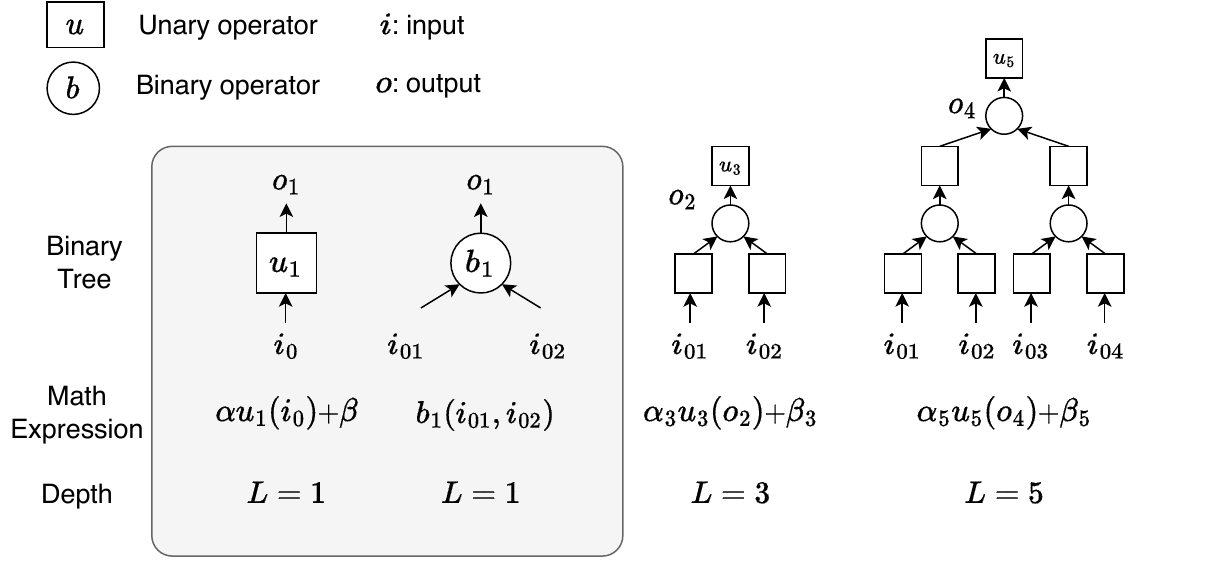}
    \caption{Computational rule of a binary tree. Each node within the binary tree holds either a unary or a binary operator. Initially, we outline the computation flow of a depth-1 tree comprising a solitary operator. Subsequently, for binary trees extending beyond a single layer, the computation process is recursively executed.}
    \label{fig:trees}
\end{figure}

\subsubsection{Finite Expressions with Binary Trees}
\label{sec:binary_tree}
FEX utilizes a binary tree $\mathcal{T}$ structure to represent finite expressions as illustrated in Figure~\ref{fig:trees}. Each tree node contains either a unary operator 
or a binary operator. The sets of unary and binary operator candidates,  $\mathbb{U}$ and $\mathbb{B}$, usually four of each kind, are selected for each problem by the user. Examples of unary and binary operators are, respectively 
$$
\sin,\exp, \log, \text{Id}, (\cdot)^2, \int\cdot\text{d} x_i, \frac{\partial\cdot}{\partial x_i}, \cdots\
\quad{\rm and}\quad 
+,-,\times,\div,\cdots.
$$
Each unary operator acts element-wise and is
equipped with scaling parameters $\alpha_i$, $i=1,\ldots,d$, and a bias parameter $\beta$ that are applied to its output element-wise, for example,
$$
\alpha_1\sin(x_1)+ \ldots +\alpha_d\sin(x_d) + \beta.
$$

The set of all parameters for all unary operators is denoted by $\bm{\theta}$.
Then, the entire expression is obtained by a preorder 
traversal of the operator sequence $\bm{e}$ of the binary tree $\mathcal{T}$. Therefore, such finite expression is denoted by $u(\bm{x};
\mathcal{T},\bm{e},\bm{\theta})$ as a function in $\bm{x}$. For a fixed $\mathcal{T}$, the maximal number of operators is
bounded from above by a constant denoted as $k_{\mathcal{T}}$. In FEX, 
$$\mathbb{S}_{\sf FEX}=\{u(\bm{x}; \mathcal{T}, \bm{e}, \bm{\theta})~|~\bm{e}\in  \mathbb{U}\cup\mathbb{B}\} 
$$ 
is the functional space in which we solve the PDE. Note that
$\mathbb{S}_{\sf FEX}\subset\mathbb{S}_{k_{\mathcal{T}}}$. The computation flow of the binary tree $\mathcal{T}$ works recursively from the leaf nodes. The unary operators at leaf nodes
are applied on the input $\bm{x}$ elementwise, and the scaling $\alpha$ transforms the dimension from $\mathbb{R}^d$ to $\mathbb{R}$. Then, the computation follows 
a bottom-up manner recursively until the flow reaches the tree root.

\subsubsection{Implementation of FEX}
Thus, the FEX solution $u(\bm{x}; \mathcal{T}, \bm{e}, \bm{\theta})$ is obtained by solving the mixed combinatorial optimization problem
of the form
\begin{align}
    \min_{\bm{e}, \bm{\theta}} \mathcal{L}(u(\cdot; \mathcal{T}, \bm{e}, \bm{\theta})).
    \label{eqn:obj}
\end{align}
To achieve this, FEX proceeds in two stages. It first optimizes the selection of operator sequence $\bm{e}$ that identifies
the structure of the true solution. Then it optimizes the parameter set $\bm{\theta}$ to minimize the funtional~\eqref{eqn:obj}.
The framework of FEX consists of four parts.
\begin{enumerate}
    \item {\it Score computation.}  To identify the structure of the solution, FEX uses a mix-order optimization algorithm to evaluate the score of the operator sequence $\bm{e}$.
    \item {\it Operator sequence generation.} FEX employs a neural network (NN) to model a controller that outputs a probability mass function to sample optimal operator sequences. 
    \item {\it Controller update.} Based on the reward feedback of generated operator sequences, the controller is updated to generate good operator sequences. 
    \item {\it Candidate optimization.} Within the search loop, FEX maintains a pool of top-performing operator sequences. After training, a fine-tuning step is performed for each candidate in the pool to obtain the best operator sequence as the approximation to the PDE solution.
\end{enumerate} 
Each of these parts will be elaborated on in the next sections.

\subsubsection{Score Computation} \label{sec:score}
The score of an operator sequence $\bm{e}$ is an essential part of training, as it
guides the controller to update parameters to output optimal probability mass functions to sample good operators. We define the score of $\bm{e}$, $S(\Be)$, by
\begin{align}
    S(\Be) := \big(1+L(\Be)\big)^{-1}, \quad{\rm where}\quad L(\Be):=\min \{\mathcal{L}(u(\cdot; \mathcal{T}, \Be, \bm{\theta}))|\bm{\theta}\}.
    \label{eqn:orgscore}
\end{align}
As $L(\Be)$ approaches zero, $S(\Be)$ increases up to 1. To efficiently evaluate the score $S(\Be)$,
 the following hybrid mix-order optimization approach is utilized for updating the parameter $\bm{\theta}$. Let $\bm{\theta^{\Be}}_0$ be the initial guess for $\bm{\theta}$ for given $\bm{e}$. First, $T_1$ steps of a first-order optimization algorithm (e.g., the stochastic gradient descent~\cite{rumelhart1986learning} or Adam~\cite{kingma2014adam}) are performed resulting in $\bm{\theta^{\Be}}_{T_1}$. Then $T_2$ steps of a second-order optimization algorithm (e.g., Newton's method~\cite{avriel2003nonlinear} or BFGS~\cite{fletcher2013practical}) are made resulting in $\bm{\theta^{\Be}}_{T_1+T_2}$. Finally, the score of 
the operator sequence $\Be$ is obtained as
\begin{align}
    S(\Be) \approx \big(1+\mathcal{L} (u(\cdot; \mathcal{T}, \Be, \bm{\theta}_{T_1+T_2}^{\Be}))\big)^{-1}. 
    \label{eqn:score}
\end{align}

\subsubsection{Operator Sequence Generation}
The goal of the controller is to output operator sequences with high scores during training. The controller $\bm{\chi}$ with parameters $\Phi$ will be denoted by $\bm{\chi}_\Phi$. For an operator sequence
$\Be$ with $s$ nodes, the controller $\bm{\chi}_\Phi$ outputs  probability mass functions $\bm{p_{\Phi}}^i$, $i = 1,\cdots,s$. Then, the operator $e_j$ is sampled from $\bm{p_{\Phi}}^j$.
In addition, the $\epsilon$-greedy strategy~\cite{sutton2018reinforcement} is used to encourage exploration in the operator set. With probability $\epsilon<1$, $e_i$ is sampled from a uniform distribution of the operator 
set and with probability $1-\epsilon$, $e_i$ is sampled from $\bm{p}_\Phi^i$.

\subsubsection{Controller Update}
The goal of the controller $\bm{\chi}_\Phi$ is to output optimal probability mass functions, from which the operator sequence $\Be$ with high scores
are highly likely to be sampled. We model the controller $\bm{\chi}_\Phi$ as a neural network parameterized by $\Phi$. The training objective is to maximize the expected score 
of a sampled operator sequence $\Be$,  i.e.
\begin{align}
    \mathcal{J}(\Phi):=\mathbb{E}_{\Be \sim \bm{\chi}_\Phi} S(\Be).
    \label{eqn:expect}
\end{align}
The derivative of \eqref{eqn:expect} with respect to $\Phi$ is 
\begin{align}
    \nabla_\Phi\mathcal{J}(\Phi)=\mathbb{E}_{\Be \sim \bm{\chi}_\Phi} \left\{S(\Be)\sum_{i=1}^s \nabla_\Phi \log\left(\bm{p}_\Phi^i(e_i)\right)\right\},
    \label{eqn:expectgrad}
\end{align}
where $\Bp_\Phi^i(e_i)$ is the probability of the sampled $e_i$. Let $N$ denote the batch size. The batch $\{\Be^{(1)}, \Be^{(2)}, \cdots, \Be^{(N)}\}$ is sampled under $\bm{\chi}_{\Phi}$ each time.
Then the expectation~\eqref{eqn:expectgrad} can be approximated by
\begin{align}
    \nabla_\Phi\mathcal{J}(\Phi)\approx \frac{1}{N}\sum_{k=1}^N \left\{S(\Be^{(k)})\sum_{i=1}^s \nabla_\Phi \log\left(\Bp_\Phi^i(e_i^{(k)})\right)\right\}.
    \label{eqn:mcmcavg}
\end{align}
In turn, the model parameter $\Phi$ is updated by gradient ascent, i.e., $\Phi \leftarrow \Phi+\eta \nabla_\Phi\mathcal{J}(\Phi)$. However, in practice, the goal is to obtain the operator sequence $\Be$ with the highest score, instead of optimizing the 
average scores of all generated operator sequences. Therefore, following~\cite{petersen2021deep} we consider
\begin{align}
    \mathcal{J}(\Phi)=\mathbb{E}_{\Be \sim \bm{\chi}_\Phi} \left\{S(\Be)|S(\Be)\geq S_{\nu, \Phi}\right\},
    \label{eqn:expectriskseeking}
\end{align}
where $S_{\nu, \Phi}$ represents the $(1-\nu)\times 100\%$-quantile of the score distribution generated by $\bm{\chi}_{\Phi}$. In a discrete form,
the gradient computation becomes
\begin{align}
    \nabla_\Phi\mathcal{J}(\Phi)\approx \frac{1}{N}\sum_{k=1}^N \left\{(S(\Be^{(k)})-\hat{S}_{\nu, \Phi})\mathbbm{1}_{\{S(\Be^{(k)})\geq \hat{S}_{\nu, \Phi}\}}\sum_{i=1}^s \nabla_\Phi \log\left(\Bp_\Phi^i(e_i^{(k)})\right)\right\},
    \label{eqn:mcmcrisk}
\end{align}
where $\mathbbm{1}$ is an indicator function that takes value $1$ if the condition is true and otherwise 0, and $\hat{S}_{\nu, \Phi}$ is the $(1-\nu)$-quantile of the scores in $\{S(\Be^{(i)})\}_{i=1}^N$.

\subsubsection{Candidate Optimization}
As introduced in Section~\ref{sec:score}, the score of $\Be$ is based on the optimization of a nonconvex function. Therefore, the score obtained by coarse-tuning with $T_1 + T_2$ iterations may not be
a good indicator of whether $\Be$ recovers the underlying structure of the solution. Therefore, it is important to keep a pool $\mathbb{P}$ of fixed size $K$, which adaptively keeps the  top $K$ candidate operator sequences $\Be$.
After the search is finished, for each $\Be \in \mathbb{P}$, the objective function $\mathcal{L}(u(\cdot; \mathcal{T}, \Be, \bm{\theta}))$ is fine-tuned over $\bm{\theta}$ using a first-order algorithm with a small learning rate for $T_3$ iterations.

\subsection{The Solution Model for the Committor Problem} \label{sec:FEX_committor}

\begin{figure}[ht]
    \centering
    \includegraphics[width=0.80\linewidth]{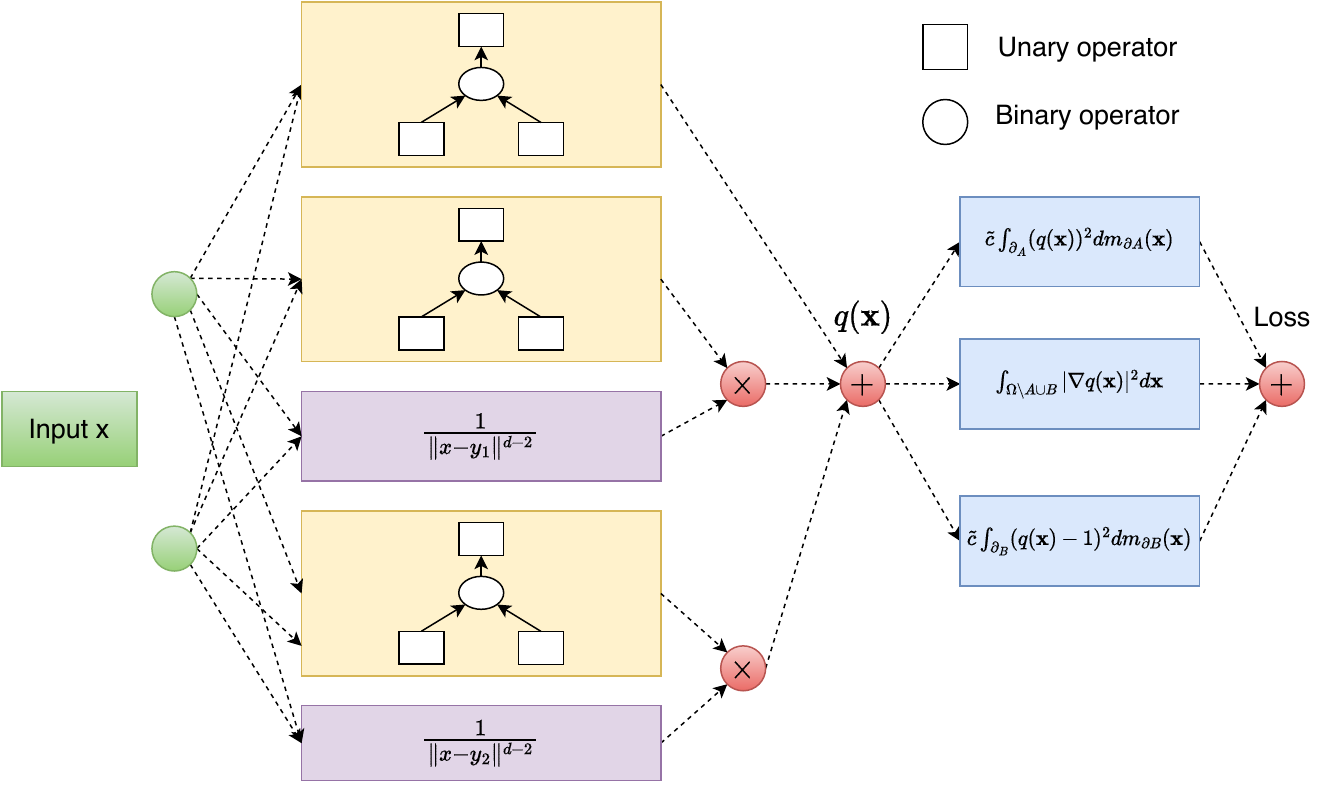}
    \caption{The representation of the computation flow of solving for~\eqref{eqn:regularized_variational_form}. The committor function $q$ is represented by the summation of three ``FEX trees'', two of which are weighted with $\frac{1}{|\bx|^{d-2}}$ type singularities. }
    \label{fig:fex_variational}
\end{figure}

In this paper, we parameterize the committor function $q(\bx)$ by a FEX binary tree
\begin{align*}
q_{\Phi} := \{q(\bm{x}; \mathcal{T}, \bm{e}, \bm{\theta})|\bm{e}, \bm{\theta}\}.
\end{align*}
We follow the setup in~\cite{khoo2019solving} and choose the variational formulation~\eqref{eqn:regularized_variational_form} as $\mathcal{L}$. As noted in~\cite{khoo2019solving}, in the high-temperature regime, i.e., when $T \rightarrow \infty$ thus $\beta \rightarrow 0$, the backward Kolmogorov equation \eqref{eqn:backward_kolmogorov} converges 
to Laplace's equation with a Dirichlet boundary condition. Therefore, the solution near boundaries $\partial A$ and $\partial B$ are dictated asymptotically by the fundamental solution
\begin{equation}
    \label{eqn:fundamental_sol}
    \Phi(\mathbf{x}):= \begin{cases}-\frac{1}{2 \pi} \log |\mathbf{x}| & (d=2), \\ \frac{\Gamma(d / 2)}{(2 \pi)^{d / 2}|\mathbf{x}|^{d-2}} & (d \geq 3).\end{cases}
\end{equation}
Considering the singular behavior of the committor function, we model the committor function as
\begin{equation}
    \label{eqn:committor_singular}
    q{(\bx)} = q_{1} S_A(\bx - \by^A) + q_{2} S_B(\bx - \by^B) + q_{3},
\end{equation}
where $\by^A$ and $\by^B$ are the centers of $A$ and $B$, $S_A(\bx - \by^A)$ and $S_B(\bx - \by^B)$ are fundamental solutions~\eqref{eqn:fundamental_sol}, and $q_1, q_2$ and $q_3$ are three FEX binary trees to be optimized. The 
computation flow is illustrated in Figure~\ref{fig:fex_variational}.

\section{Numerical Experiments}
\label{sec:numericalexp}
In this section, the performance of FEX on a collection of benchmark test problems is examined. All test problems are set up in high-dimensional ambient spaces. All these problems admit variable changes reducing them to low-dimensional problems. We demonstrate that FEX is effective at identifying this low-dimensional structure automatically.
Thus, our main objective in this section is to illustrate the capabilities of FEX in the following ways.


\begin{enumerate}
\item FEX demonstrates comparable or higher accuracy compared to the neural network method. 
\item FEX excels in identifying the low-dimensional structure inherent in each problem. 
\item Once FEX successfully identifies the low-dimensional structure, we can achieve arbitrary accuracy by solving the reduced low-dimensional problem~\eqref{eqn:backward_kolmogorov} using spectral methods or finite element methods. 
\end{enumerate}

The committor problem \eqref{eqn:backward_kolmogorov} seldom admits an analytical solution. Therefore, in order to evaluate the accuracy of FEX one needs to find a highly accurate solution by another method. We use the finite element method (FEM) for this purpose.

Prior to delving into the subsequent benchmark problems, we shall provide a concise overview of the convergence analysis associated with the finite element method (FEM). The significance of discussing the error estimate for the finite element method holds a two-fold importance in our numerical experiments. Firstly, in line with our previously stated primary objective, the comparison of the relative error between the neural network method and the ground truth solution, as well as the relative error between FEX and the ground truth solution, requires the use of a highly accurate finite element method as the reference solution. Secondly, as part of our third objective, once the low-dimensional structure of the problem~\eqref{eqn:backward_kolmogorov} has been identified, we propose to use the spectral method and finite element method to solve the low-dimensional ODE (or PDE)~\eqref{eqn:backward_kolmogorov}, both of which possess a robust theoretical foundation for achieving arbitrary accuracy.

As mentioned in Theorem 5.4 of~\cite{larsson2003partial}, the finite element method utilizing piecewise linear basis functions to address elliptic PDEs within a convex two-dimensional domain also showcases quadratic convergence of the numerical solution towards the exact solution with respect to step size $h$.

It is important to note that Theorem 5.4 presented in~\cite{larsson2003partial} provides preliminary error estimates with certain limitations. Firstly, these theorems were derived in the context of a simplified elliptic PDE with only the second order term. Therefore, its applicability to our problem~\eqref{eqn:backward_kolmogorov} needs to be considered carefully. Additionally, Theorem 5.4~\cite{larsson2003partial} assumes a convex domain $\Omega$, whereas our benchmark problems involve non-convex domains. Despite these limitations, these theorems serve as a starting point for determining an appropriate step size $h$ and assessing the overall accuracy of our benchmark results. In the subsequent benchmark problems, where the finite element method is employed as the reference solution, we specifically utilize piecewise linear basis functions, but the choice of the step size $h$ is tailored for each individual case.

As a data-driven solver, FEX requires a set of training points as input. The variational loss functional \eqref{eqn:regularized_variational_form} assumes two kinds of training points: $N_{\sf bdry}$ boundary points lying on $\partial A\cup\partial B$ and $N_{\Omega_{AB}}$ interior points lying in $\Omega_{AB}$.
In our experiments, we use $N_{\sf bdry} = 2000$ and keep the ratio 
$$
\frac{N_{\sf bdry}}{N_{\Omega_{AB}}}=:2\alpha
$$
between $1/10$ to $1/100$. 
To evaluate the accuracy, we use the relative $L_2$ error
\begin{equation*}
    E=\frac{\left\|q_\theta-q\right\|_{L_2(\mu)}}{\|q\|_{L_2(\mu)}},
\end{equation*}
where $q_\theta$ represents the numerical solution obtained through either the neural network (NN) or FEX methods, while $q$ denotes the reference solution characterized by its high accuracy.

The loss functional \eqref{eqn:regularized_variational_form} can be rewritten as a single expectation
\begin{equation}
    \label{eqn:loss_expectation}
    \mathcal{L}(q) = \mathbb{E}_{\nu} \left(|\nabla q(\bx)|^2\chi_{\Omega_{AB}} (\bx) + \tilde{c}q(\bx)^2 \chi_{\partial A}(\bx) + \tilde{c}(q(\bx)-1 )^2\chi_{\partial B}(\bx)\right),
\end{equation}
where $\nu$ is the mixture measure $    \nu (\bx) = \rho(\bx) +m_{\partial A}(\bx) + m_{\partial B}(\bx)$ with $\rho$ being the invariant density given by \eqref{eqn:inv_pdf} and $m_{\partial A}$, $m_{\partial B}$ are uniform measures on the boundaries $\partial A$ and $\partial B$ respectively.

For FEX, we employ a binary tree of depth either $L = 3$ or $L = 5$. From Figure~\ref{fig:trees}, it can be observed that a depth-3 tree has two leaf nodes, whereas a depth-5 tree has four leaf nodes. These leaf nodes are labeled in a left-to-right order. For example, in the case of a depth-5 tree, the leftmost leaf node is labeled as leaf 1, and the rightmost leaf node is labeled as leaf 4.
We select the binary set $\mathbbm{B}=\{+,-,\times\}$ and the unary set $\mathbbm{U}=\{0, 1, \text{Id}, (\cdot)^2, (\cdot)^3, (\cdot)^4, \exp, \sin, \cos, \tanh, \text{sigmoid}\}$ to form the mathematical equation. Notably, we include the $\tanh$ and sigmoid functions in $\mathbbm{U}$ due to the often observed sharp transitions in committor functions within high-dimensional spaces. Additionally, we include the $(\cdot)^2$ operator, as the committor function is often related to the spherical radius of the problem. We present a selection of algebraic formulas identified by FEX in the following examples, while the remaining equations can be found in the supplementary material.

\subsection{The Double-Well Potential with Hyperplane Boundary}
The first example features the committor problem \eqref{eqn:backward_kolmogorov} with the double-well potential ~\cite{khoo2019solving,chen2023committor}
\begin{equation}
\label{eqn:pw_potential}
    V(\mathbf{x})=\left(x_1^2-1\right)^2+0.3 \sum_{i=2}^d x_i^2,
\end{equation}
and 
\begin{equation*}
\label{eqn:doublewell_boundary}
    A=\left\{\bx \in \R^d \mid x_1 \leq-1\right\}, \quad B=\left\{\bx \in \R^d \mid x_1 \geq 1\right\}.
\end{equation*}
We set $d=10$.
Note that this problem is effectively one-dimensional, as the committor depends only on the first component of $\bx$, $x_1$, and is the solution to 
\begin{equation}
\label{eqn:pw_sol}
    \frac{d^2 q\left(x_1\right)}{d x_1^2}-4\beta  x_1\left(x_1^2-1\right) \frac{d q\left(x_1\right)}{d x_1}=0, \quad q(-1)=0, \quad q(1)=1.
\end{equation}
The solution to \eqref{eqn:pw_sol} is given by
\begin{equation*}
    \label{eqn:exact1}
    q(x_1) = \frac{\int_{-1}^{x_1} e^{\beta (y^2-1)^2}dy}{\int_{-1}^{1} e^{\beta (y^2-1)^2}dy}.
\end{equation*}

Now we let FEX find out that the solution depends only on $x_1$. We model the committor function $q(\bx)$ as a single depth-3 FEX tree $\mathcal{J}(\bx)$ as there is no singularity. When the temperature $\beta^{-1}$ is low, sampling from the invariant density gives too few samples near the transition state at $\bx = 0$. So, following~\cite{khoo2019solving}, we sample $x_1$ uniformly on $[-1,1]$, and $(x_2,\cdots,x_d)$ from a $(d-1)$-dimensional Gaussian distribution. This sampling density is accounted for in the Monte Carlo integration of the loss functional \eqref{eqn:loss_expectation} resulting in 
\begin{align}
     \mathcal{L}(q)& = \frac{1}{N_{\Omega_{AB}}}\sum_{j=1}^{N_{\Omega_{AB}}} 
    \frac{\left|\nabla q(\mathbf{x}_j)\right|^2 \exp \left(-\beta\left((x_1)_j^2-1\right)^2\right)}{\int_{-1}^1 \exp \left(-\beta\left((x_1)_j^2-1\right)^2\right) d x_1} \notag \\
    &+ \frac{\tilde{c}}{N_{\sf bdry}}\sum_{j=1}^{N_{\sf bdry}/2} \left[q(\mathbf{y}_j)^2 + 
     (1-q(\mathbf{z}_j))^2\right],
     \label{eqn:modified_interior_loss}
\end{align}
where $\bx_j\in\Omega_{AB}$, $j=1,\ldots,N_{\Omega_{AB}}$, and $\mathbf{y}_j \in\partial A$, $\mathbf{z}_j\in\partial B$.

\begin{figure}[ht]
    \centering
    \begin{subfigure}[b]{0.45\textwidth}
      \centering
      \includegraphics[width=\textwidth]{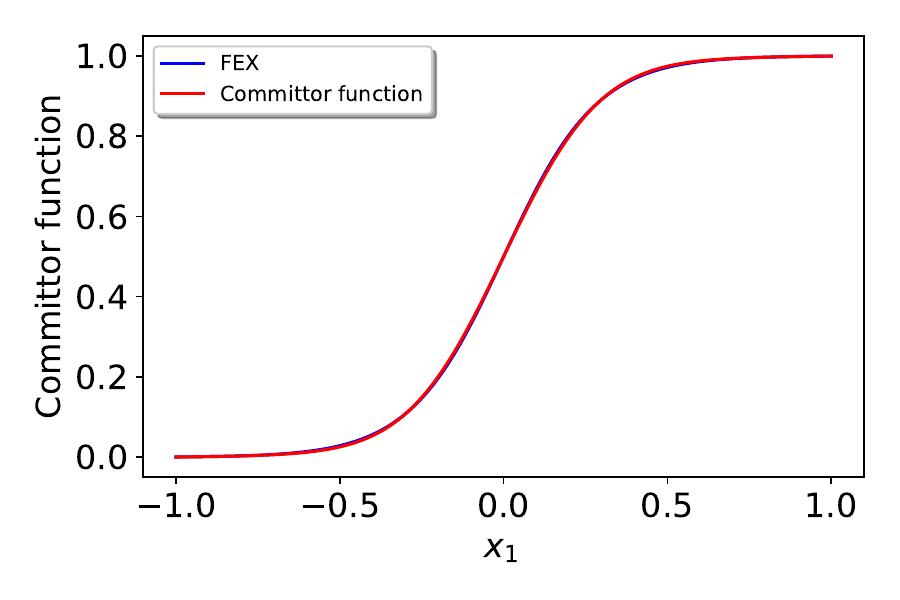}
      \caption{Committor function and FEX}
      \label{fig:pw_FEX_T_0_2}
    \end{subfigure}
    \hfill
    \begin{subfigure}[b]{0.45\textwidth}
      \centering
      \includegraphics[width=\textwidth]{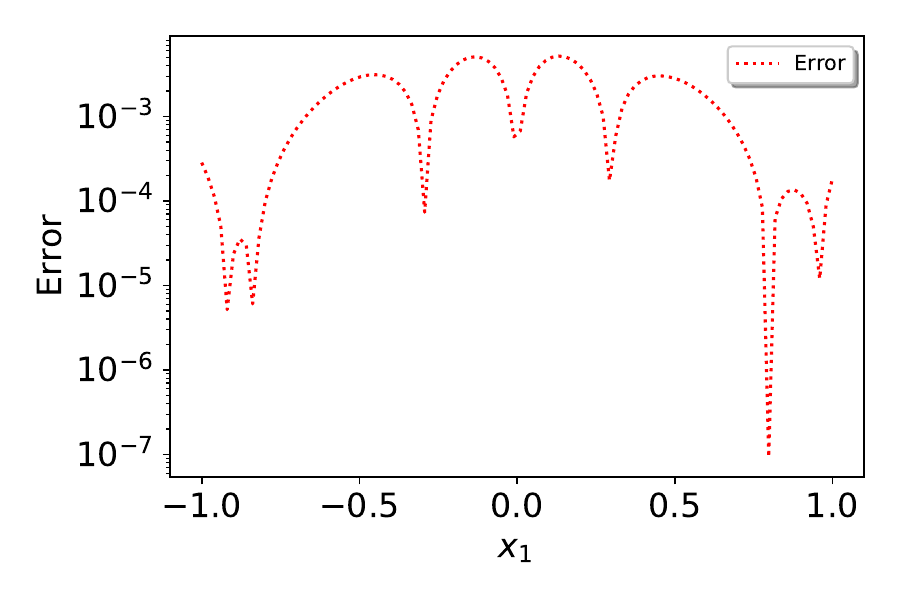}
      \caption{Error of FEX when solving \eqref{eqn:regularized_variational_form}}
      \label{fig:pw_FEX_error_T_0_2}
    \end{subfigure}
  
  
    \caption{The committor function for the double-well potential along $x_1$ dimension when $\beta^{-1}=0.2$ for an arbitrary $(x_2,\cdots,x_d)$ with $d=10$.}
    \label{fig:pw_T_0.2}
  \end{figure}

We consider $\beta^{-1} = 0.2$ and $\beta^{-1} = 0.05$, respectively. At $\beta^{-1} = 0.2$, the expression of $\mathcal{J}(\bx)$ found by FEX is 
\begin{align*}
    & \text{leaf 1: \texttt{Id}} \rightarrow \alpha_{1,1} x_1 + \ldots + \alpha_{1,{10}} x_{10} + \beta_1 \\
    & \text{leaf 2: \texttt{tanh}} \rightarrow \alpha_{2,1} \tanh(x_1) + \ldots + \alpha_{2,{10}} \tanh(x_{10}) + \beta_2 \\
    & \mathcal{J}(\bx) = \alpha_{3} \tanh(\text{leaf 1} + \text{leaf 2}) + \beta_3,
\end{align*}
where $\alpha_{3} = 0.5,~ \beta_3 = 0.5$. The rest parameters in leaf 1 and leaf 2 are summarized in Table~\ref{tab:2well_coef_tree_T_0.2}. It is evident that only the coefficients for $x_1$ and $\tanh(x_1)$ are nonzero. The formula for $\beta^{-1} = 0.05$  computed by FEX is placed in the supplementary material. We plot the FEX committor function and its error with the true committor function in Figure~\ref{fig:pw_T_0.2}. We also summarize the numerical results in Table~\ref{tab:doublewell_T_0.2} and Table~\ref{tab:doublewell_T_0.05}.

\begin{table}
    \center
    \begin{tabular}{lccccccccccc}
    \Xhline{3\arrayrulewidth} 
    \text{node} &  $\alpha_1$ & $\alpha_2$ & $\alpha_3$ & $\alpha_4$ & $\alpha_5$ & $\alpha_6$ & $\alpha_7$ &$\alpha_8$ &$\alpha_9$ &$\alpha_{10}$ &$\beta$\\
    \hline
    leaf 1: Id   & $\mathbf{1.6798}$ & $0.0$ & $0.0$ & $0.0$ & $0.0$ & $0.0$ & $0.0$ &$0.0$ &$0.0$ &$0.0$ &$0.0$\\
    \hline
    leaf 2: \texttt{tanh}   & $\mathbf{1.9039}$ & $0.0$ & $0.0$ & $0.0$ & $0.0$ & $0.0$ & $0.0$ &$0.0$ &$0.0$ &$0.0$ &$0.0$ \\
    \Xhline{3\arrayrulewidth}
    \end{tabular}
    \caption{Coefficients of leaves of depth-3 FEX binary tree $\mathcal{J}$ for double-well potential with hyperplane boundary problem when $T = 0.2$. We would like to emphasize that, in the fine-tuning stage of the FEX algorithm, we employ coefficient filtering mentioned in Algorithm~\ref{alg:workflow} with a threshold of $\tau = 0.05$ to effectively eliminate trivial coefficients.
    Therefore, the result clearly shows that FEX identifies the structure of the committor function $q(x_1)$, enabling post-processing techniques, such as spectral method.}
    \label{tab:2well_coef_tree_T_0.2}
\end{table}




Clearly, FEX successfully identifies that the committor $q(\bx)$ depends only on the first coordinate $x_1$. This allows us to reduce the problem to 1D and 
obtain a machine precision solution  to \eqref{eqn:pw_sol} using the Chebyshev spectral method (\cite{trefethen2000spectral}, program  {\tt p13.m}). Since equation \eqref{eqn:pw_sol} is linear, the boundary conditions are enforced by decomposing the solution as $q(x_1) = u(x_1) + u_b(x_1)$ where $u(x_1)$ satisfies the same equation with homogeneous boundary conditions and $u_b(x_1)$ is any smooth user-chosen function that satisfies $u_b(-1) =0$ and $u_b(1) = 1$, e.g., $u_b(x_1) = 0.5(x_1+1)$. Then the problem for $u(x_1)$ becomes $\mathcal{L}_1u = -\mathcal{L}_1u_b$ where $\mathcal{L}_1$ is the differential operator in the left-hand side of \eqref{eqn:pw_sol}. About $80$ collocation points are enough to achieve the error of the order of $10^{-16}$. The computed solution in the form of a Chebyshev sum is evaluated at any point $x_i\in[-1,1]$ using Clenshaw's method \cite{GilSeguraTemme} (see Section 3.7.1).

\begin{table}
\centering
    \begin{tabular}{ccccccc}
    \Xhline{3\arrayrulewidth} 
    \text{method} & E & $\tilde{c}$ & \begin{tabular}{c}
    No. of samples in \\
    $\Omega_{AB}$
    \end{tabular} & $\alpha$ & \begin{tabular}{c}
    No. of testing samples
    \end{tabular} \\
    \hline
    \text{NN}~\cite{khoo2019solving} & $5.40 \times 10^{-3}$ & 50 & $2.0 \times 10^{4}$ & 1/20 & $1.0 \times 10^5$ \\
    \text{FEX} & $3.51 \times 10^{-3}$  & 50 & $2.0 \times 10^{4}$ & 1/20 & $1.0 \times 10^5$ \\
    \Xhline{3\arrayrulewidth}
    \end{tabular}
    \caption{Results for the double-well potential  with hyperplane boundary when $\beta^{-1} = 0.2$.}
    \label{tab:doublewell_T_0.2}
\end{table}
\begin{table}
\centering
    \begin{tabular}{ccccccc}
    \Xhline{3\arrayrulewidth} 
    \text{method} & E & $\tilde{c}$ & \begin{tabular}{c}
    No. of samples in \\
    $\Omega_{AB}$
    \end{tabular} & $\alpha$ & \begin{tabular}{c}
    No. of testing samples
    \end{tabular} \\
    \hline
    \text{NN}~\cite{khoo2019solving} & $1.20 \times 10^{-2}$ & 0.5 & $2.0 \times 10^{4}$ & 1/20 & $1.0 \times 10^5$ \\
    \text{FEX} & $5.50 \times 10^{-3}$  & 0.5 & $2.0 \times 10^{4}$ & 1/20 & $1.0 \times 10^5$ \\
    \Xhline{3\arrayrulewidth}
    \end{tabular}
    \caption{Results for the double-well potential with hyperplane boundary when $\beta^{-1} = 0.05$.}
    \label{tab:doublewell_T_0.05}
\end{table}

\subsection{The Double-Well Potential with Sublevel Set Boundary}
In this example, we still consider the double-well potential~\eqref{eqn:pw_potential}, but with a more challenging boundary condition. Namely, we consider the sublevel sets boundaries,
\begin{equation*}
    \label{eqn:levelset_boundary}
        A=\left\{\bx \in \R^d \mid V(\bx) < V_a, ~x_1 < 0\right\}, \quad B=\left\{\bx \in \R^d \mid V(\bx) < V_b,~x_1 > 0\right\}.
\end{equation*}
In this case, the committor function $q$ solves the following equation
\begin{align*}
    \label{eqn:levelset_truth}
    &\nabla \cdot \left( \exp(-\beta V)\cdot \nabla q \right) + \frac{d-2}{r} \exp(-\beta V) \frac{dq}{dr} = 0, \\
    &\left.q(\bx)\right|_{\partial A}=0,\quad\left.q(\bx)\right|_{\partial B}=1,\notag
\end{align*}
where the committor function depends on $x_1$ and $r = \sqrt{x_2^2 + \cdots + x_d^2}$, i.e., $q = q(x_1, r)$. The values of the parameters are as follows:
$V_a = 0.2,~V_b = 0.2,~d = 10$. Given the non-convex nature of the problem domain, the finite element method (FEM) is employed as the benchmark solution. We choose the mesh size $h = 0.02$. The discrepancy error between the finite element solution and the exact solution in the $\|L\|_{\infty}$-norm is estimated to be approximately of the order of $O(10^{-4})$ to $O(10^{-5})$.
When $\beta$ is small, the committor $q$ is more heavily dependent on the radius $r$, whereas as $\beta$ increases, it effectively depends only on $x_1$. Therefore, we experiment with different values of $\beta$ and compare the performance of FEX and that of a neural network-based solver. We use a depth-5 FEX tree. More details are provided in the supplementary material. 

\begin{table}[ht]
\centering
    \begin{tabular}{lcccccc}
    \Xhline{3\arrayrulewidth} 
    $\beta$ & \text{method} & E & $\tilde{c}$ & \begin{tabular}{c}
    No. of samples in \\
    $\Omega_{AB}$
    \end{tabular} & $\alpha$ & \begin{tabular}{c}
    No. of testing samples
    \end{tabular} \\
    \hline
    3.0 & \text{NN} & $0.372$ & 50 & $2.0 \times 10^{4}$ & 1/20 & $2.0 \times 10^{4}$ \\
    3.0 & \text{FEX} & $0.358$  & 50 & $2.0 \times 10^{4}$ & 1/20 & $2.0 \times 10^{4}$ \\
    10.0 & \text{NN} & $0.021$ & 50 & $2.0 \times 10^{4}$ & 1/20 & $2.0 \times 10^{4}$ \\
    10.0 & \text{FEX} & $0.009$  & 50 & $2.0 \times 10^{4}$ & 1/20 & $2.0 \times 10^{4}$ \\
    \Xhline{3\arrayrulewidth}
    \end{tabular}
    \caption{Results for the double-well potential with sublevel sets boundaries.}
    \label{tab:2well_levelset}
\end{table}

\begin{figure}[ht]
    \centering
    \begin{subfigure}[b]{0.45\textwidth}
      \centering
      \includegraphics[width=\textwidth]{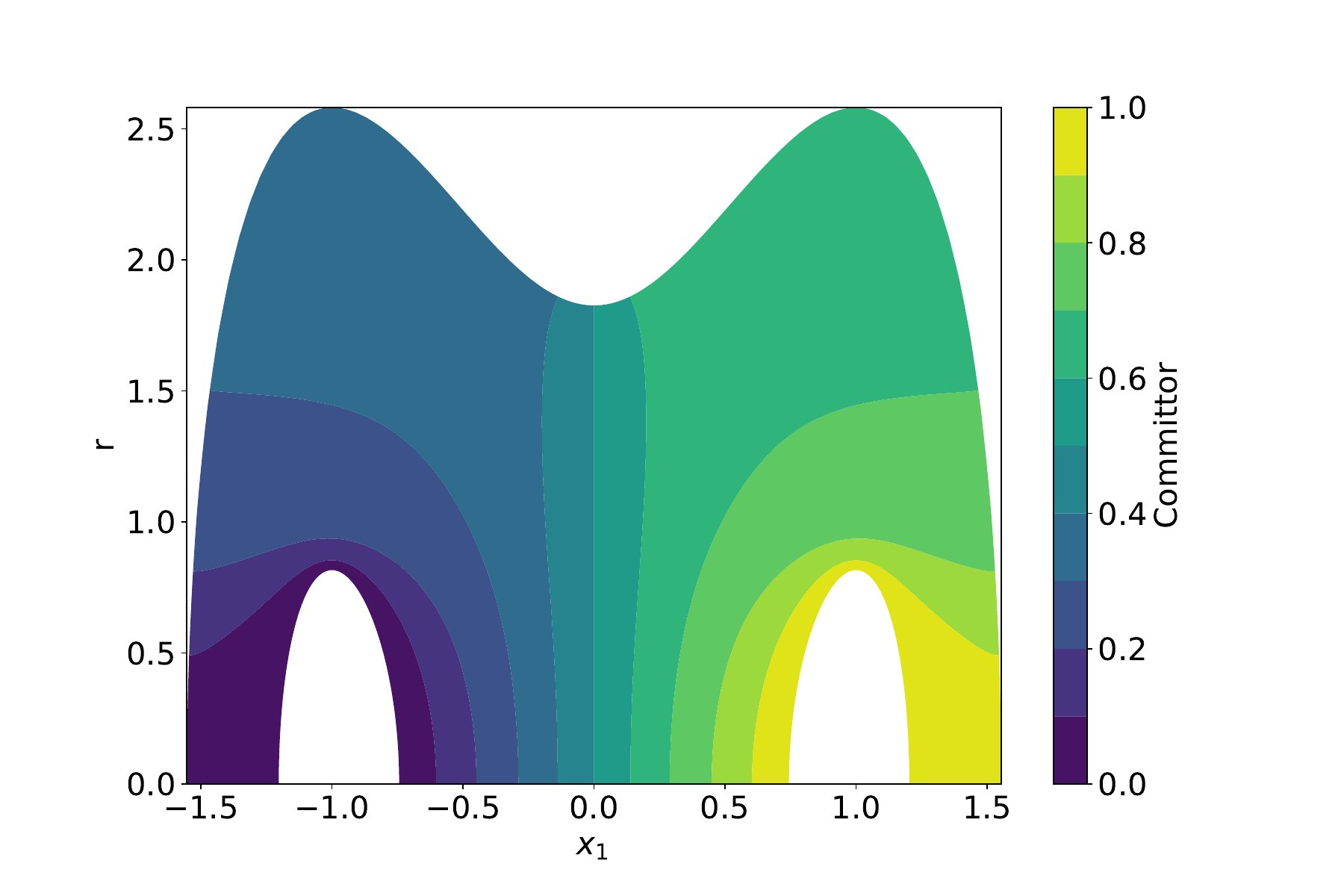}
      \caption{$\beta = 3.0$ committor function (FEM)}
      \label{fig:fem_levelset_beta_3}
    \end{subfigure}
    \hfill
    \begin{subfigure}[b]{0.45\textwidth}
      \centering
      \includegraphics[width=\textwidth]{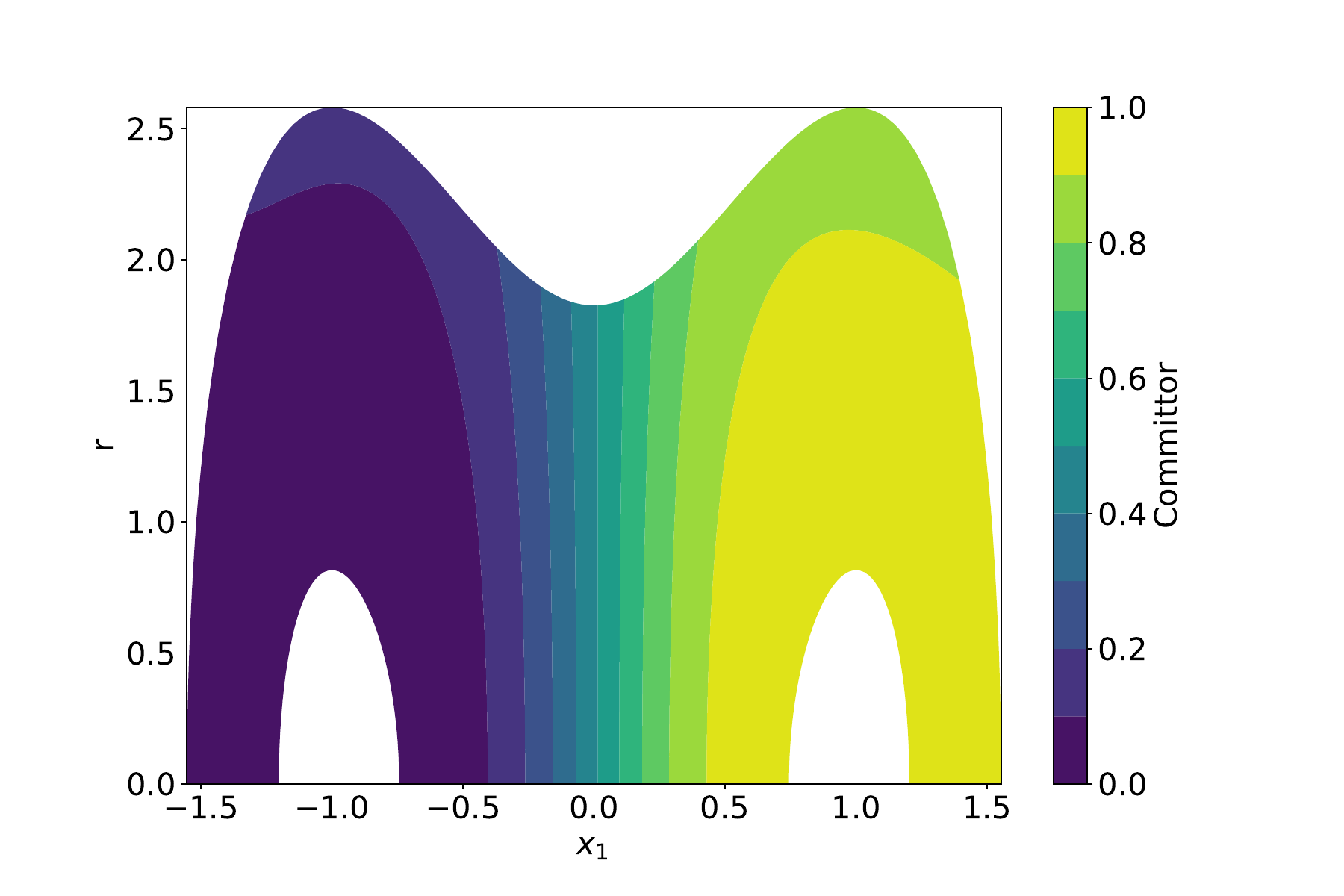}
      \caption{$\beta = 3.0$ committor function (FEX)}
      \label{fig:fex_levelset_beta_3}
    \end{subfigure}
  
    \medskip
  
    \begin{subfigure}[b]{0.45\textwidth}
      \centering
      \includegraphics[width=\textwidth]{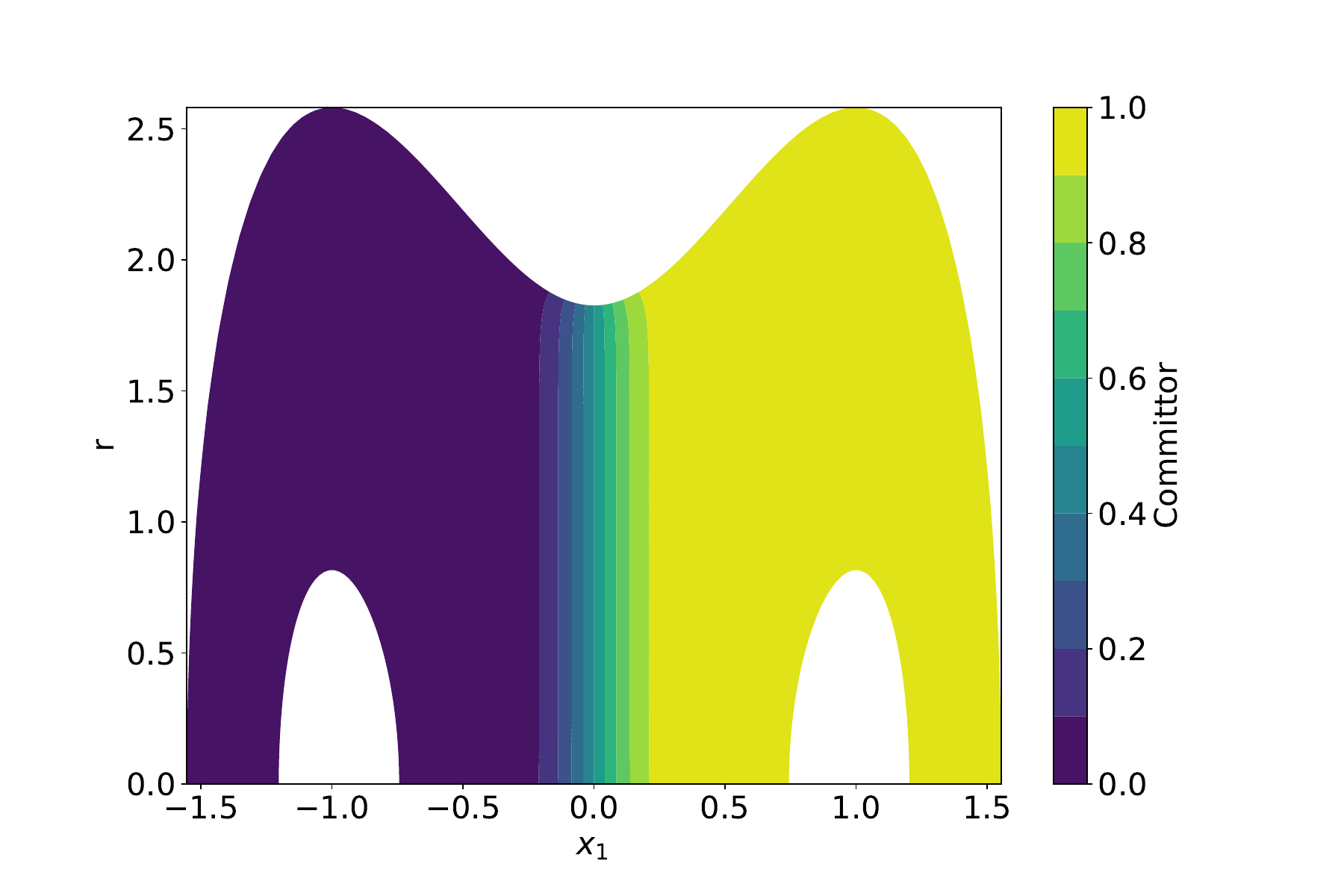}
      \caption{$\beta = 10.0$ committor function (FEM)}
      \label{fig:fem_levelset_beta_10}
    \end{subfigure}
    \hfill
    \begin{subfigure}[b]{0.45\textwidth}
      \centering
      \includegraphics[width=\textwidth]{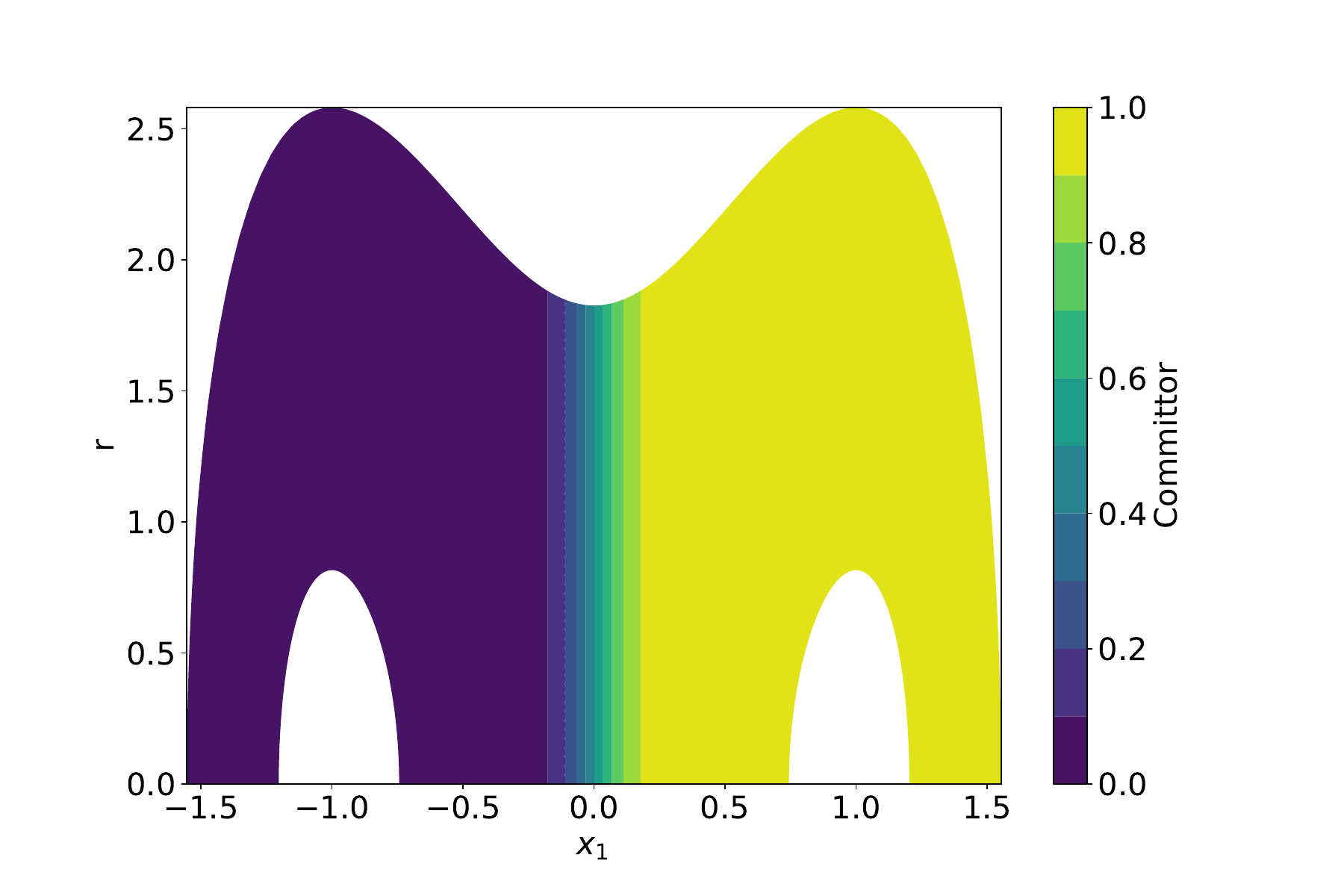}
      \caption{$\beta = 10.0$ committor function (FEX)}
      \label{fig:fex_levelset_beta_10}
    \end{subfigure}
    \caption{Committor function for the double-well potential with sublevel sets boundary. As $\beta$ increases, the problem transforms from a 2D problem to a 1D problem, and FEX can capture such behavior of the committor function.}
    \label{fig:sublevel_set}
  \end{figure}

  \begin{figure}[ht]
    \centering
    \begin{subfigure}[b]{0.45\textwidth}
      \centering
      \includegraphics[width=\textwidth]{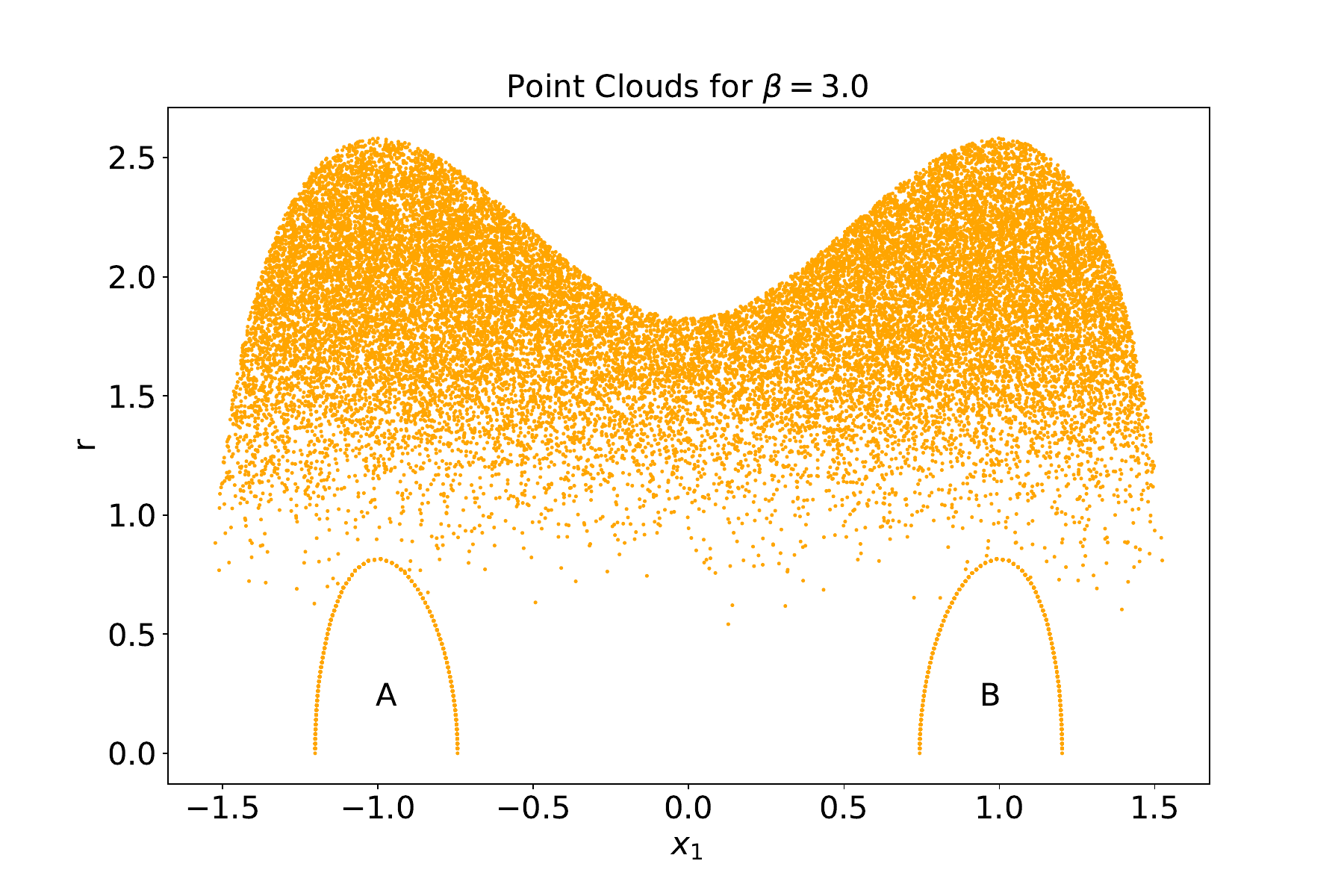}
      \caption{}
    \end{subfigure}
    \hfill
    \begin{subfigure}[b]{0.45\textwidth}
      \centering
      \includegraphics[width=\textwidth]{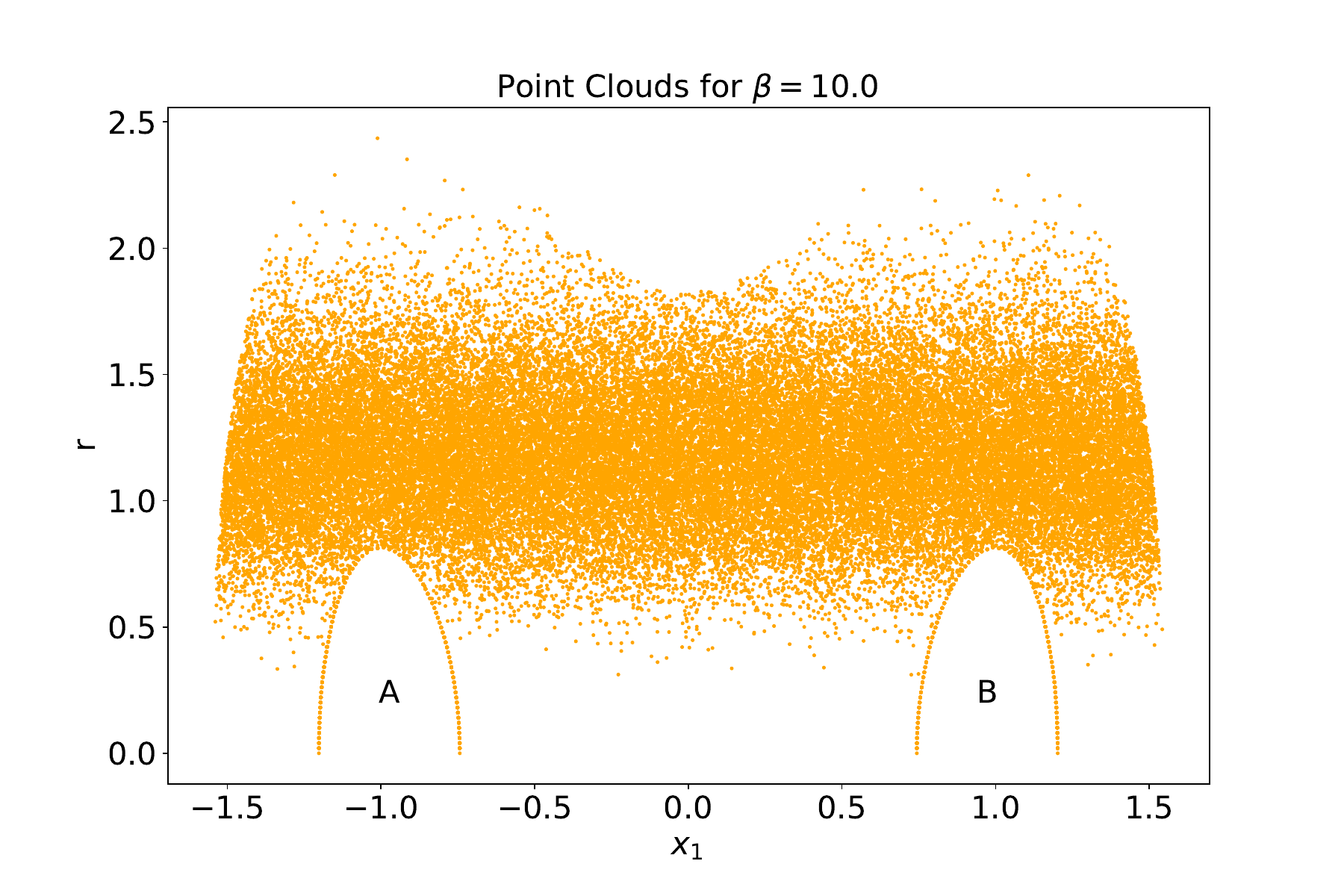}
      \caption{}
    \end{subfigure}
  
    \caption{Point clouds for double-well potential with sublevel set boundaries. Left: $\beta = 3.0$. Right: $\beta = 10.0$.} 
    \label{fig:point_sublevel}
  \end{figure}

From Table \ref{tab:2well_levelset}, we observe that FEX consistently outperforms NN for various temperatures. Furthermore, when $\beta = 3.0$, FEX could indeed identify the structure
of the committor function, i.e., $q(x_1, r)$. While when $\beta = 10.0$, the problem essentially reduces to a 1D problem, and FEX is able to identify such change of pattern. In comparison, NN based methods are unable 
to retrieve such geometric information. The numerical solutions by FEM and by FEX are displayed in Fig. \ref{fig:sublevel_set}. From Figure.~\ref{fig:point_sublevel}, we observe that for $\beta = 3.0$, point clouds generated from the invariant density~\eqref{eqn:inv_pdf} are sparse near and between boundaries $A$ and $B$. In contrast, for $\beta = 10.0$, the data coverage is more sufficient. This is challenging for any data-driven approaches, which explains the relatively large errors for both FEX and neural network-based methods for $\beta = 3.0$. Nevertheless, FEX is still able to identify the low-dimensional structure of the problem in such a data-scarce regime, i.e., $q:= q(x_1, r)$.

\subsection{Concentric Spheres}
In this example, we consider the committor function describing the transition process between a pair of concentric spheres, with the potential
\begin{equation}
    \label{eqn:2sphere_potential}
    V(\bx) = 10|\bx|^2,
\end{equation}
and the regions 
\begin{equation*}
    \label{eqn:2sphere_region}
    A=\left\{\bx \in \mathbb{R}^d\mid |\bx| \geq a\right\}, \quad B=\left\{\bx \in \mathbb{R}^d\mid |\bx| \leq b\right\}.
\end{equation*}

\begin{figure}[ht]
    \centering
    \begin{subfigure}[b]{0.45\textwidth}
      \centering
      \includegraphics[width=\textwidth]{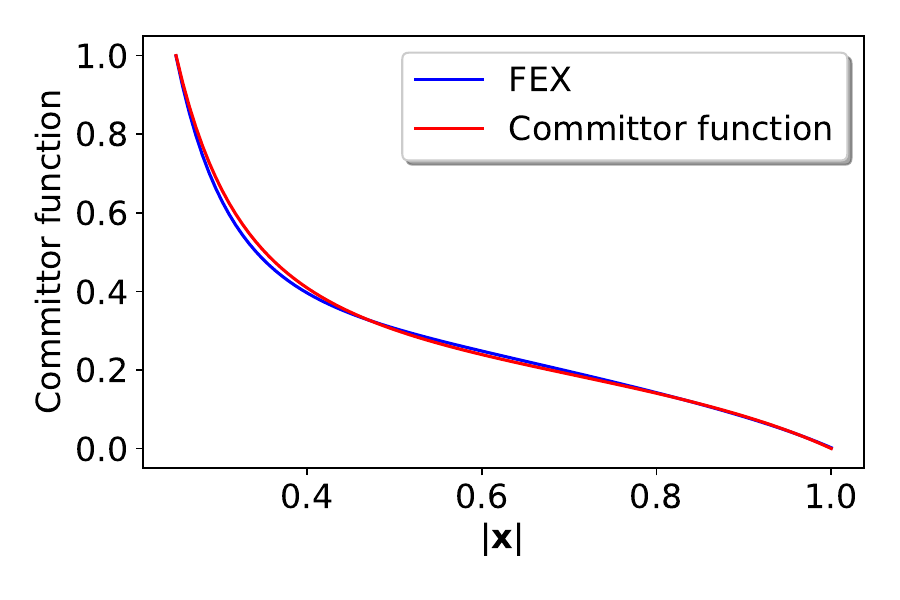}
      \caption{Committor function and FEX}
      \label{fig:2sphere_FEX_T_0.2}
    \end{subfigure}
    \hfill
    \begin{subfigure}[b]{0.45\textwidth}
      \centering
      \includegraphics[width=\textwidth]{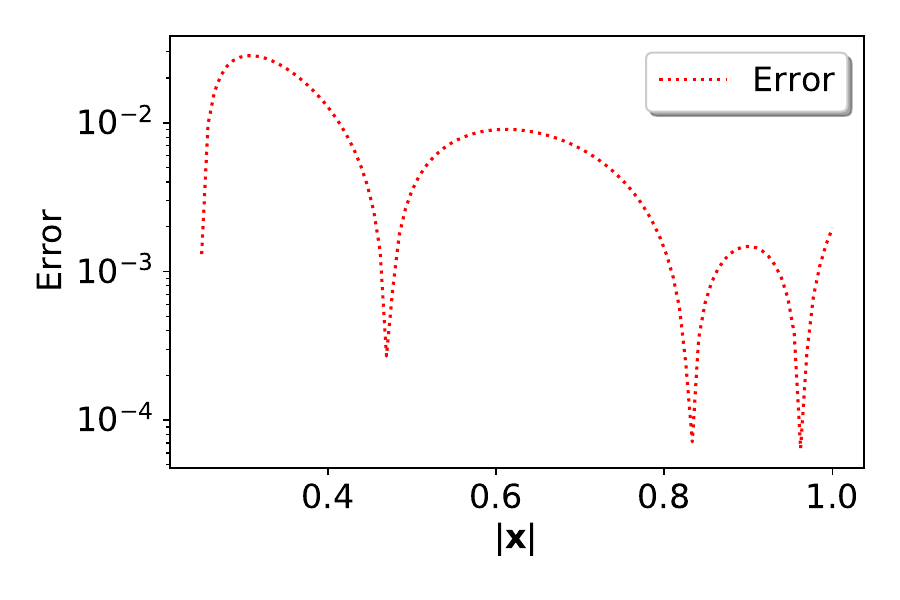}
      \caption{Error of FEX when solving \eqref{eqn:regularized_variational_form}}
      \label{fig:2sphere_FEX_error_T_0.2}
    \end{subfigure}
      \caption{The committor function for the concentric spheres as a function of $|\bx|$.}
      \label{fig:2sphere_committor}
  

  \end{figure}

Since the equilibrium distribution is proportional to Gaussian distribution, we can readily obtain samples from the normal distribution. 
We obtain the data on the boundaries $\partial A, \partial B$ by sampling from Gaussian distribution and rescale to have norm $a$ or $b$. 
The values of parameters are: $T = 2, d = 6, a = 1, b = 0.25$.
The true solution $q(\bx) := q(r)$ satisfies the following ODE:
\begin{align}
    &\label{eqn:true_sol_2sphere}
    \frac{d^2q(r)}{dr^2} + \frac{d-1}{r}\frac{dq(r)}{dr} - \beta \frac{dq}{dr}\frac{dV}{dr} = 0,\\
    &\left.q(r)\right|_{r=a}=0,\quad\left.q(r)\right|_{r=b}=1,\notag
\end{align}
where $r = |\bx|$ is the radius, $\beta = 1/T$ is the inverse of temperature. Similarly to ~\eqref{eqn:pw_sol}, ~\eqref{eqn:true_sol_2sphere} is solvable.
In this example, the committor function displays a singular behavior, i.e., $q \sim 1/|\bx|^{d-2}, d \geq 3$. Therefore, we parameterize the committor function $q(\bx)$ as 
\begin{equation*}
    q{(\bx)} = \mathcal{J}_1(\bx) \cdot \frac{1}{|\bx|^{d-2}} + \mathcal{J}_2(\bx),
\end{equation*}
where $\mathcal{J}_1(\bx)$ and $\mathcal{J}_2(\bx)$ are two depth-3 FEX binary trees. We leave the FEX formula in the supplementary material. The formula reduces to
\begin{equation*}
    q{(\bx)} = q(r) := \frac{0.0020}{r^{0.5d-1}} + 0.6016(0.6054 - 0.5800r^2)-0.0340,
\end{equation*}
Therefore, FEX can successfully identify that the committor function $q(\bx)$ only depends on the spherical radius $r$.
This allows us to solve \eqref{eqn:true_sol_2sphere}  using the Chebyshev spectral method once again as described in the previous section and obtain the numerical solution with the machine precision. The numerical results are provided in Table~\ref{tab:2sphere}.
In addition, we plot the committor function as a function of $|\bm{x}|$ in Figure~\ref{fig:2sphere_committor}.

\begin{table}[ht]
\centering
    \begin{tabular}{ccccccc}
    \Xhline{3\arrayrulewidth} 
    \text{method} & E & $\tilde{c}$ & \begin{tabular}{c}
    No. of samples in \\
    $\Omega_{AB}$
    \end{tabular} & $\alpha$ & \begin{tabular}{c}
    No. of testing samples
    \end{tabular} \\
    \hline 
    \text{NN}~\cite{khoo2019solving} & $5.30 \times 10^{-2}$ & 530 & $3.0 \times 10^{4}$ & 1/30 & $1.0 \times 10^5$ \\
    \text{FEX} & $3.20 \times 10^{-2}$  & 530 & $3.0 \times 10^{4}$ & 1/30 & $1.0 \times 10^5$ \\
    \Xhline{3\arrayrulewidth}
    \end{tabular}
    \caption{Results for the concentric spheres example.}
    \label{tab:2sphere}
\end{table}

\subsection{Rugged Mueller's Potential}
In this example, we consider the committor function in the rugged Mueller's potential
\begin{equation*}
    \label{eqn:rugged_muller}
    V(\bx)=\tilde{V}\left(x_1, x_2\right)+\frac{1}{2 \sigma^2} \sum_{i=3}^d x_i^2,
\end{equation*}
where
\begin{equation*}
    \label{eqn:2d_rugged_muller}
    \tilde{V}\left(x_1, x_2\right)=\sum_{i=1}^4 D_i e^{a_i\left(x_1-X_i\right)^2+b_i\left(x_1-X_i\right)\left(x_2-Y_i\right)+c_i\left(x_2-Y_i\right)^2}+\gamma \sin \left(2 k \pi x_1\right) \sin \left(2 k \pi x_2\right),
\end{equation*}
where $\gamma = 9$ and $k = 5$ determine the ruggedness of the 2-dimensional rugged Mueller's potential $\tilde{V}\left(x_1, x_2\right)$, $\sigma = 0.05$ controls the extent of the harmonic potential in dimensions $x_3, \cdots, x_d$, and $d = 10$. All other parameters are listed below, which is consistent with \cite{khoo2019solving}.
\begin{align*}
{\left[a_1, a_2, a_3, a_4\right] } & =[-1,-1,-6.5,0.7], & {\left[b_1, b_2, b_3, b_4\right] } & =[0,0,11,0.6], \\
{\left[c_1, c_2, c_3, c_4\right] } & =[-10,-10,-6.5,0.7], & {\left[D_1, D_2, D_3, D_4\right] } & =[-200,-100,-170,15], \\
{\left[X_1, X_2, X_3, X_4\right] } & =[1,0,-0.5,-1], & {\left[Y_1, Y_2, Y_3, Y_4\right] } & =[0,0.5,1.5,1].  & &
\end{align*}
We focus on the domain $\Omega = [-1.5,1] \times[-0.5,2] \times \mathbb{R}^{d-2}$ and the regions $A$ and $B$ are two cylinders:

\begin{align*}
    \label{domain_rm}
    & A=\left\{\bx \in \mathbb{R}^d \mid \sqrt{\left(x_1+0.57\right)^2+\left(x_2-1.43\right)^2} \leq 0.3\right\} \\
    & B=\left\{\bx \in \mathbb{R}^d \mid \sqrt{\left(x_1-0.56\right)^2+\left(x_2-0.044\right)^2} \leq 0.3\right\}.
\end{align*}

The ground truth solution is obtained by the finite element method, which is obtained by solving \eqref{eqn:backward_kolmogorov} on uniform grid in 2 dimensions with the potential $\tilde{V}$, the domain $\tilde{\Omega} = [-1.5, 1] \times [-0.5, 2]$, and the region $\tilde{A}$, $\tilde{B}$ projected by $A$ and $B$ on the $x_1 x_2$-plane, respectively. The chosen step size $h$ for the finite element method is on the order of $O(10^{-3})$, resulting in a discretization error on the order of $O(10^{-6})$ to $O(10^{-7})$ when compared to the ground truth solution. Notably, these error magnitudes are considerably smaller than the numerical errors incurred by the numerical solvers NN and FEX. Hence, we designate it as a reference solution of notable precision and accuracy. In this case, there are singularities
present in regions $A$ and $B$, so we parameterize the committor function as
\begin{equation}
    q{(\bx)} = \mathcal{J}_1 \log((x_1 + 0.57)^2 + (x_2 - 1.43)^2) + \mathcal{J}_2 \log((x_1 - 0.56)^2 + (x_2 - 0.044)^2) + \mathcal{J}_0,
\end{equation}
where $\mathcal{J}_1, \mathcal{J}_2, \mathcal{J}_0$ are three depth-5 FEX binary trees, respectively.

\begin{figure}[ht]
    \centering
    \begin{subfigure}[b]{0.45\textwidth}
      \centering
      \includegraphics[width=\textwidth]{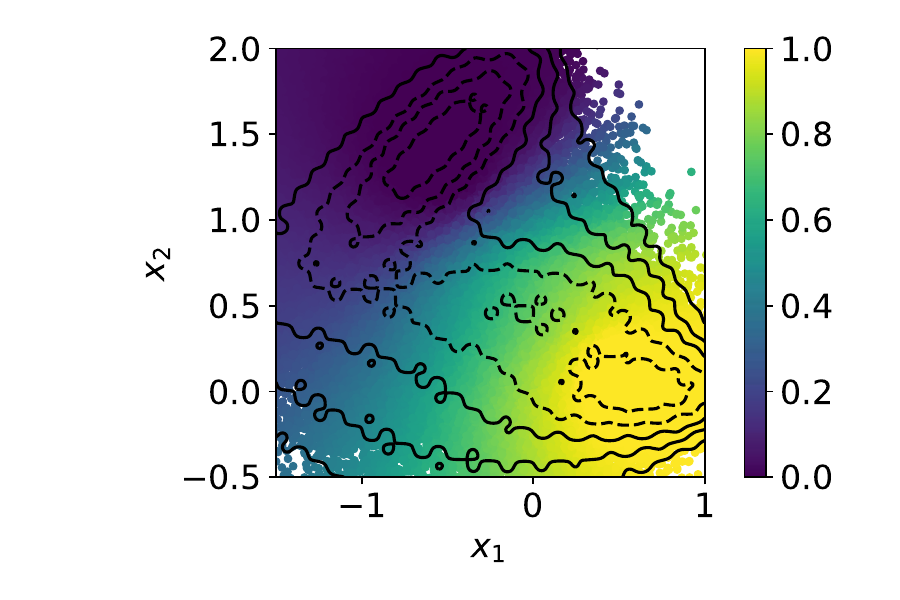}
      \caption{$T = 40$ committor function (FEM)}
      \label{fig:RM_T_40_True}
    \end{subfigure}
    \hfill
    \begin{subfigure}[b]{0.45\textwidth}
      \centering
      \includegraphics[width=\textwidth]{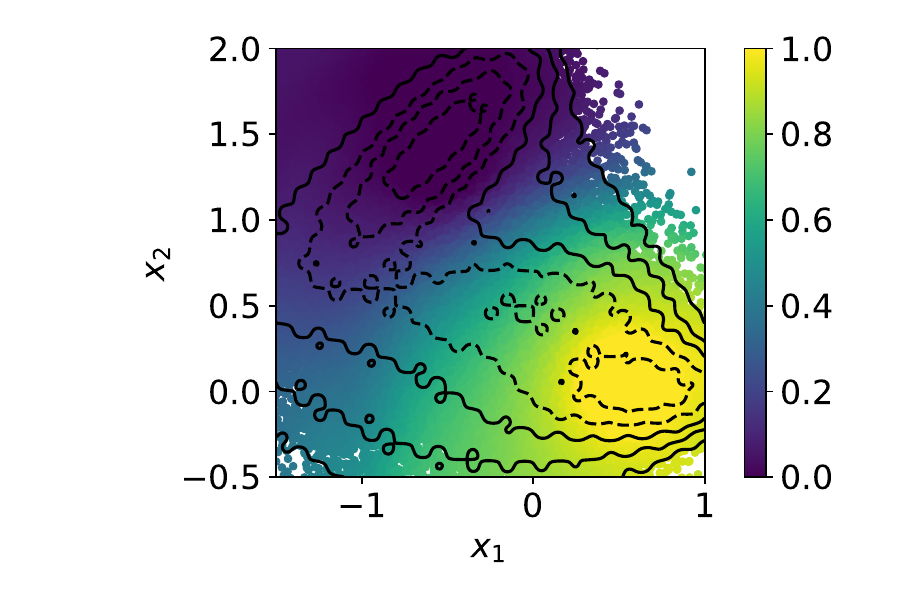}
      \caption{$T = 40$ committor function (FEX)}
      \label{fig:RM_T_40_FEM}
    \end{subfigure}
  
    \medskip
  
    \begin{subfigure}[b]{0.45\textwidth}
      \centering
      \includegraphics[width=\textwidth]{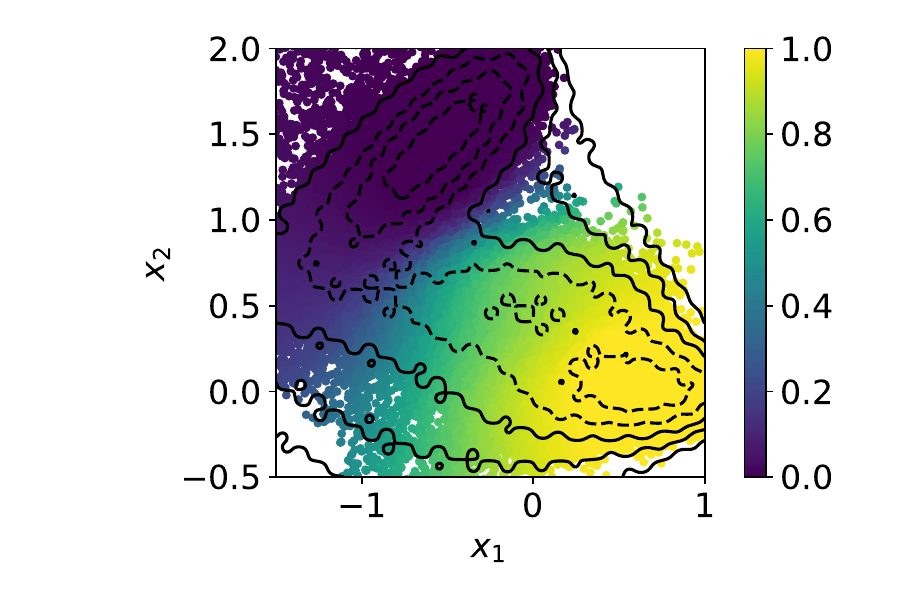}
      \caption{$T = 22$ committor function (FEM)}
      \label{fig:RM_T_22_True}
    \end{subfigure}
    \hfill
    \begin{subfigure}[b]{0.45\textwidth}
      \centering
      \includegraphics[width=\textwidth]{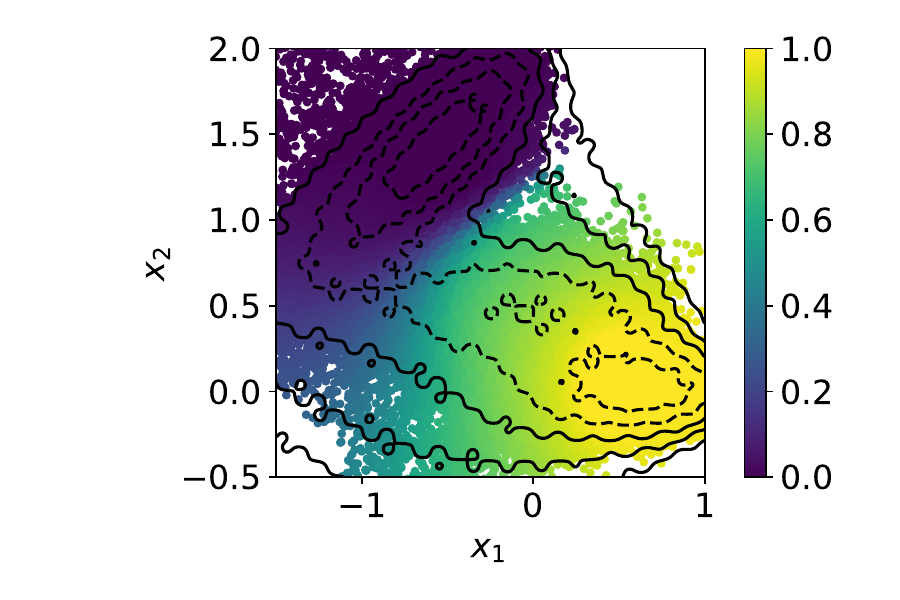}
      \caption{$T = 22$ committor function (FEX)}
      \label{fig:RM_T_22_FEM}
    \end{subfigure}
    \caption{Committor function for the rugged Mueller's potential on a 2-dimensional plane, with level curves of potential $\tilde{V}$.}
    \label{fig:RM}
  \end{figure}

We consider $T = 40$ and $T = 22$, and delay the FEX formula in the supplementary material, and these formulas show that FEX identifies the solution only varies with the first two coordinates $x_1, x_2$. Therefore, we further simplify the problem and instead solve~\eqref{eqn:backward_kolmogorov} with $V(\bx)=\tilde{V}\left(x_1, x_2\right)$.
We utilize finite element method to obtain the committor function $q(\bx) = q(x_1,x_2)$. We plot the committor function $q(x_1,x_2)$ on the $x_1x_2$ plane in Figure~\ref{fig:RM}. We also summarize the numerical results in Table~\ref{tab:muller}.

  \newpage
  \subsection{Butane}
 
\begin{figure}[ht]
    \centering
    \begin{subfigure}[b]{0.45\textwidth}
      \centering
      \includegraphics[width=\textwidth]{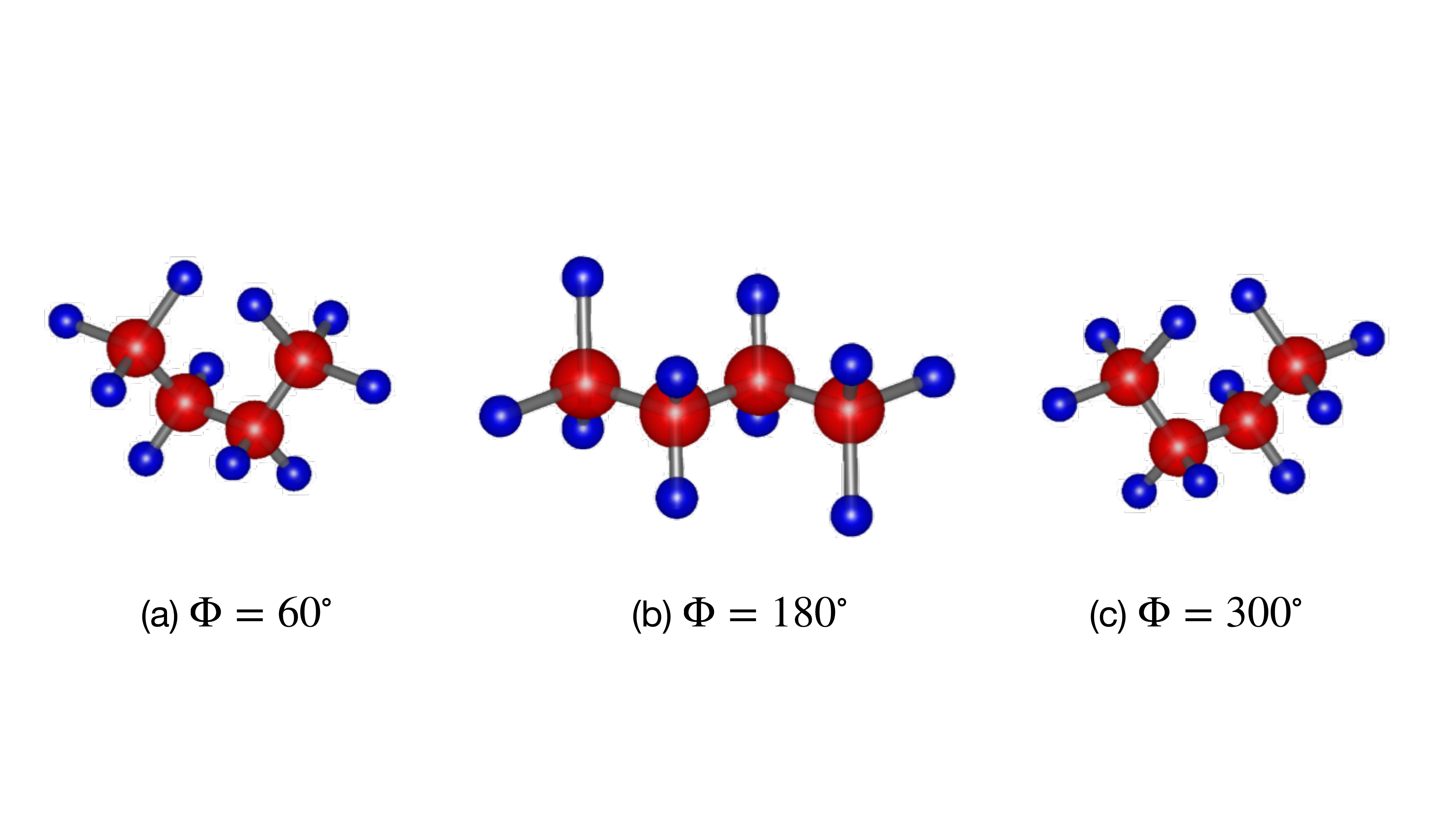}
      \caption{}
      \label{fig:butane_config}
    \end{subfigure}
    \hfill
    \begin{subfigure}[b]{0.45\textwidth}
      \centering
      \includegraphics[width=\textwidth]{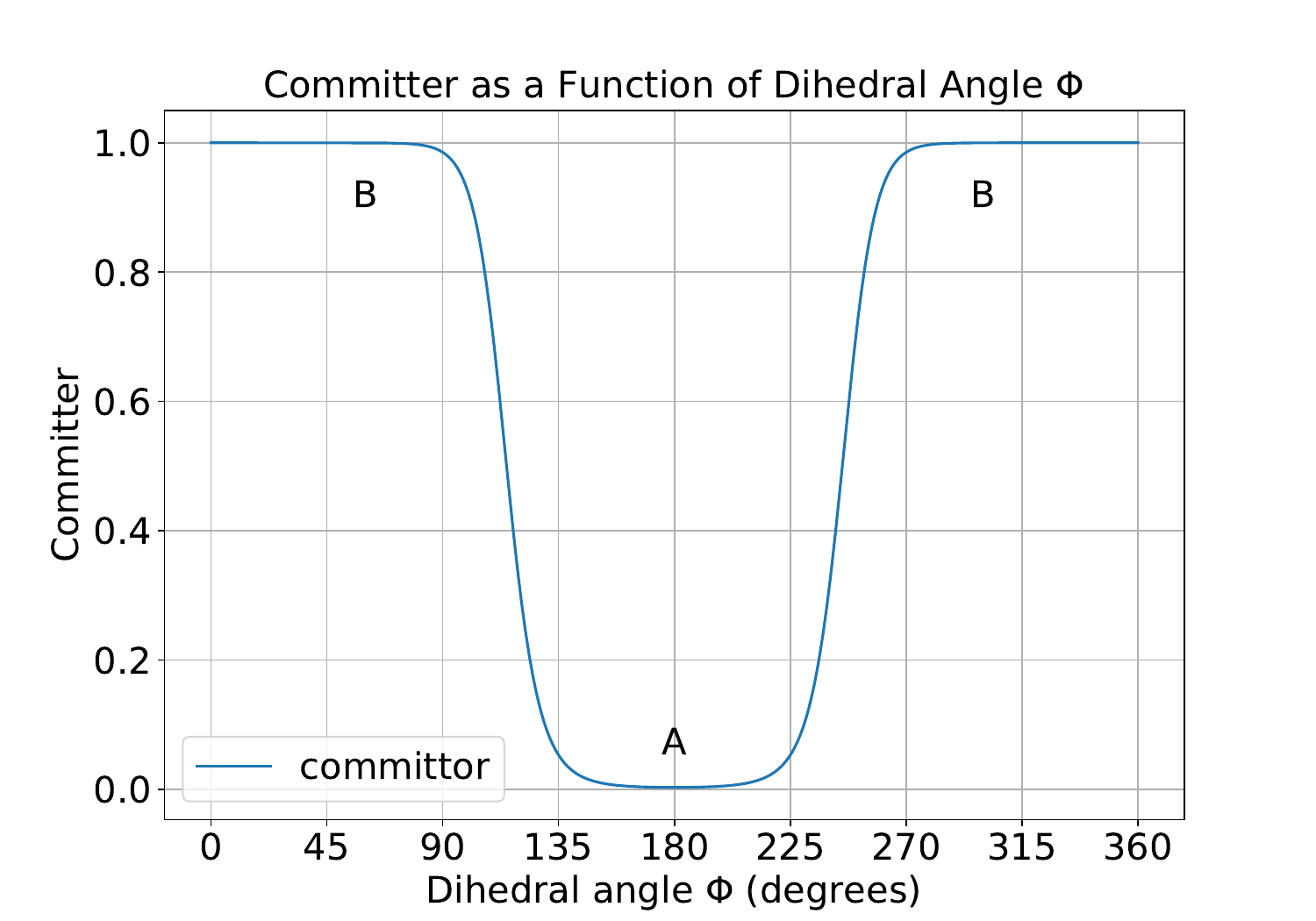}
      \caption{}
      \label{fig:butane_committor}
    \end{subfigure}
  
    \caption{Left: Butane configuration at different dihedral angles $\Phi$. Right: Committor as a function of dihedral angle $\Phi$. }
    \label{fig:butane}
  \end{figure}

  Butane is a small hydrocarbon molecule with chemical formula C$_4$H$_{10}$. A chain-like butane molecule  is a popular test problem for quantifying transition times between its three metastable states shown in Fig. \ref{fig:butane}(a) \cite{Zuckerman2001TransitionEI,Klus2015OnTN,banisch2020}. Hydrogen atoms are light compared to carbons and are typically ignored in any type of coarse-graining. The carbons are connected by covalent bonds resulting in nearly constant distances between the bonded atoms. The main mode of the motion of butane is due to the changes in the dihedral angle, i.e. the angle between the two planes passing through the first, second, and third carbons and the second, third, and fourth carbons along the carbon backbone of butane. The metastable states are typically defined as intervals surrounding the values of the dihedral angle of 60, 180, and 300 degrees lifted into the coordinate space of carbons. The covalent (or interbond) angles in change in the process of thermal motion but significantly less as the dihedral angle.

 Our goal in this example to test is FEX is able to identify that the committor function depends on the dihedral angle and is independent of the covalent angles.

  The input data for the butane were obtained from a molecular dynamics simulation at the temperature of 300K using OpenMM~\cite{eastman2017openmm}. The simulation
  consisted of $10^7$ time steps. The timestep was set to 2 femtoseconds. The trajectory snapshots data consisted of $N$ snapshots, with each snapshot having a dimension of $D = 12$ (4 atoms $\times$ 3 coordinates for $x,y,z$.). Subsampling was performed every 100 steps, resulting in data points at
  0.2 picoseconds intervals.

  The diheral angle and the covalent angles can be expressed via the vectors $C_i C_{i+1}$ defined as
  \begin{equation*}
      C_i C_{i+1} = \left(x_{i+1} - x_i, y_{i+1} - y_i, z_{i+1} - z_{i}\right) \quad \text{for} \quad i = 1, 2, 3.
  \end{equation*}
 Then the cosine of the dihedral angle $\Phi$ is equal to
  \begin{equation*}
      \cos{\Phi} = \frac{\mathbf{n}_1 \cdot \mathbf{n}_2}{\|\mathbf{n}_1\|_2 \|\mathbf{n}_2\|_2},
  \end{equation*}
  where $\mathbf{n}_1, \mathbf{n}_2$ are the normal vectors of planes, spanned by $C_1 C_2$ and $C_2 C_3$ , and $C_2 C_3$ and $C_3 C_4$ respectively, i.e.,
  \begin{equation*}
      \begin{aligned}
          \mathbf{n}_1 = C_1 C_2 \times C_2 C_3, \\
          \mathbf{n}_2 = C_2 C_3 \times C_3 C_4.
      \end{aligned}
  \end{equation*}
  The covalent angles $\theta_1$ and $\theta_2$ are defined as
  \begin{equation*}
      \begin{aligned}
          \cos{\theta_1} = \frac{C_1C_2 \cdot C_2C_3}{\|C_1C_2\|_2 \|C_2C_3 \|_2}, \\
          \cos{\theta_2} = \frac{C_2C_3 \cdot C_3C_4}{\|C_2 C_3\|_2 \|C_3 C_4 \|_2}.
      \end{aligned}
  \end{equation*}

  The reactant state $A$ is defined as 1000 point clouds with dihedral angle $\Phi = 180\pm 0.5730^\circ
  $. The product state $B$ is defined as the union of 1000 point clouds with dihedral angle  $\Phi = 60\pm 0.5730^\circ
  $ and  $\Phi = 300\pm 0.5730^\circ$. 
We hypothesize that the committor depends on the cosines of the diherdal and covalent angles, i.e., $q:= q(\cos{\Phi}, \cos{\theta_1}, \cos{\theta_2})$.

  The configuration of butane molecule at the reactant and product states are visualized in Fig.~\ref{fig:butane_config}. The expression of $\mathcal{J}(\bx)$ found by FEX is
  \begin{align*}
      & \text{leaf 1: \texttt{$(\cdot)$}} \rightarrow 10.2078\cdot \cos{\Phi}+0.0000\cdot\cos{\theta_1}+0.0000\cdot \cos{\theta_2}+4.2090 \\
      & \text{leaf 2: \texttt{$(\cdot)$}} \rightarrow 0.0000\cdot \cos{\Phi} +0.0000\cdot \cos{\theta_1}+0.0000\cdot\cos{\theta_2}+ 0.0000 \\
      & \mathcal{J}(\bx) = \texttt{sigmoid}(\text{leaf 1} + \text{leaf 2}),
  \end{align*}
  Therefore, FEX is capable of identifying the low-dimensional structure of the problem that the committor depends solely on
  the dihedral angle $\Phi$. A visualization of committor as a function of the collective variable $\Phi$ is provided in Fig.~\ref{fig:butane_committor}.

\subsection{Summary}
In summary, our investigations have demonstrated that the FEX method exhibits similar or even superior accuracy compared to the neural network method across all benchmark problems. Moreover, FEX effectively captures the low-dimensional structure in each case, enabling direct and highly accurate solutions to the backward Kolmogorov equation~\eqref{eqn:backward_kolmogorov} without relying on Monte Carlo integration in the variational formulation approach~\eqref{eqn:variational_form}. Consequently, we propose employing the spectral method or finite element method to solve~\eqref{eqn:backward_kolmogorov} due to their well-established theoretical convergence rates towards the ground truth solution, enabling the attainment of arbitrary levels of accuracy.

\begin{table}[ht]
    \resizebox{\textwidth}{!}{
    \begin{tabular}{cccccccc}
    \Xhline{3\arrayrulewidth} 
    T &\text{Method} & E & $\tilde{c}$ & $h$ & \begin{tabular}{c}
    No. of samples in \\
    $\Omega_{AB}$
    \end{tabular} & $\alpha$ & \begin{tabular}{c}
    No. of testing samples
    \end{tabular} \\
    \hline
    40 & \text{NN}~\cite{khoo2019solving} &  $5.70 \times 10^{-2}$  & $3.8 \times 10^{2}$ & 0.005 & $7.4 \times 10^4$ & 1 / 74 & $7.4 \times 10^4$ \\
    40 & \text{FEX} & $5.01 \times 10^{-2}$ & $3.8 \times 10^{2}$ & 0.005  & $7.4 \times 10^4$ & 1 / 74 & $7.4 \times 10^4$ \\
    22 & \text{NN}~\cite{khoo2019solving} & $3.70 \times 10^{-2}$   & $1.3 \times 10^{2}$ & 0.005 & $1.0 \times 10^5$ & 1 / 150 & $1.5 \times 10^5$ \\
    22 & \text{FEX} & $2.90 \times 10^{-2}$ & $1.3 \times 10^{2}$ & 0.005 & $1.0 \times 10^5$ & 1 / 150 & $1.5 \times 10^5$ \\

    \Xhline{3\arrayrulewidth}
    \end{tabular}}
    \caption{Results for the rugged Mueller's potential. Since FEX identifies the committor depends only on
    $x_1$, $x_2$, we can further use finite element method to improve the accuracy.}
    \label{tab:muller}
\end{table}


\section{Conclusion}
\label{sec:summary}
In this work, we have investigated the novel finite expression method  as a solver for high-dimensional committor functions. Our numerical results show that FEX
can achieve superior performance compared to NN solvers. Most importantly, FEX is capable of identifying the low-dimensional algebraic structure of the problem that can be used for reformulating the committer problem as a low-dimensional one and finding a highly accurate solution to it by a suitable traditional technique. 

\section*{Acknowledgement}
We thank Dr. Luke Evans for providing us with the data for butane.
H. Y. was partially supported by the US National Science Foundation under awards DMS-2244988, DMS-2206333, the Office of Naval Research Award N00014-23-1-2007, and the DARPA D24AP00325-00. M.C. was partially supported by AFOSR MURI grant  FA9550-20-1-0397. Approved for public release; distribution is unlimited.

\bibliography{bibliography} 
\bibliographystyle{plain}


\end{document}